\documentclass[reqno]{amsart}

\usepackage[utf8]{inputenc}
\usepackage{times}
\usepackage[UKenglish]{babel}
\usepackage[T1]{fontenc}
\usepackage{verbatim,calligra,mathrsfs}
\usepackage[normalem]{ulem}
\usepackage{pdflscape}

\usepackage[a4paper,top=2cm,bottom=2cm,left=2cm,right=2cm,marginparwidth=1.75cm]{geometry}

\usepackage{amsmath,amssymb,amsfonts,amsthm,upgreek}
\usepackage{mathrsfs, mathtools}
\usepackage{graphicx}
\usepackage{booktabs}
\usepackage{tikz-cd}
\usepackage[colorinlistoftodos,prependcaption, textsize=small, textwidth=2.5cm]{todonotes}
\usepackage[all]{xy}
\usepackage{cancel}
\usepackage{appendix}
\usepackage{enumerate}
\usepackage{multicol}
\usetikzlibrary{calc}

\usepackage[colorlinks=true, allcolors=blue,backref=page]{hyperref}

\newcommand{\cH}{\mathcal{H}}

\newcommand{\cJ}{\mathcal{J}}

\newcommand{\cV}{\mathcal{V}}

\newcommand{\fg}{{\mathfrak g}}

\newcommand{\fh}{{\mathfrak h}}

\newcommand{\fk}{{\mathfrak k}}

\newcommand{\fm}{{\mathfrak m}}

\newcommand{\fp}{{\mathfrak p}}
\newcommand{\fq}{{\mathfrak q}}

\newcommand{\N}{\mathbb{N}}

\newcommand{\R}{\mathbb{R}}
\newcommand{\C}{\mathbb{C}}
\renewcommand{\H}{\mathbb{H}}

\newcommand{\bbS}{\mathbb{S}}

\newcommand{\gl}{\mathfrak{gl}}

\newcommand{\fsp}{\mathfrak{sp}}

\newcommand{\SO}{{\rm SO}}
\newcommand{\Sp}{{\rm Sp}}
\newcommand{\Spin}{{\rm Spin}}
\newcommand{\SU}{{\rm SU}}
\newcommand{\GL}{\mathrm{GL}}

\newcommand{\U}{{\rm U}}

\newcommand{\Cas}{\mathrm{Cas}}

\newcommand{\bangle}[1]{\langle #1 \rangle}

\newcommand{\Ad}{\mathrm{Ad}}

\renewcommand{\det}{\mathop\mathrm{det}\nolimits}

\renewcommand{\epsilon}{\varepsilon}

\newcommand{\Hom}{{\mathrm{Hom}}}
\DeclareMathOperator{\chom}{\mathscr{H}\text{\kern -3pt {\calligra\large om}}\,}

\newcommand{\Ric}{{\rm Ric}}

\newcommand{\ad}{\mathrm{ad}}

\newcommand{\diag}{\mathrm{diag}}

\newcommand{\id}{\mathrm{Id}}

\renewcommand{\Im}{\mathop{\mathrm{Im}}}

\renewcommand{\Re}{\mathop{\mathrm{Re}}}

\newcommand{\tr}{\mathop{\mathrm{tr}}\nolimits}

\newcommand{\Sym}{\mathrm{Sym}}

\renewcommand{\geq}{\geqslant}
\renewcommand{\leq}{\leqslant}

\def\<{\mathopen{}\left<}
\def\>{\right>\mathclose{}}
\def\({\mathopen{}\left(}
\def\){\right)\mathclose{}}

\usepackage{multicol, color}

\definecolor{gold}{rgb}{0.85,.66,0}
\definecolor{cherry}{rgb}{0.9,.1,.2}
\definecolor{burgundy}{rgb}{0.8,.2,.2}
\definecolor{orangered}{rgb}{0.85,.3,0}
\definecolor{orange}{rgb}{0.85,.4,0}
\definecolor{olive}{rgb}{.45,.4,0}
\definecolor{lime}{rgb}{.6,.9,0}
\definecolor{green}{rgb}{.2,.7,0}
\definecolor{grey}{rgb}{.4,.4,.2}
\definecolor{brown}{rgb}{.4,.3,.1}

\newtheorem{theorem}{Theorem}[section]
\newtheorem*{theorem*}{Theorem}
\newtheorem{proposition}{Proposition}[section]
\newtheorem{corollary}[proposition]{Corollary}
\newtheorem{lemma}[proposition]{Lemma}

\theoremstyle{remark}
\newtheorem{remark}[equation]{Remark}

\theoremstyle{definition}
\newtheorem{definition}[proposition]{Definition}

\numberwithin{equation}{section}

\title{On the Morse index of some Riemannian homogeneous fibrations}
\author[L. Lara]{Luiz Lara}
\address{Instituto de Matemática, Estatística e Computação Científica
(IMECC), Universidade Estadual de Campinas (Unicamp),
Rua Sérgio Buarque de Holanda, 651, 13083-859 Campinas, SP, Brazil}
\address{Laboratoire de Mathématiques de Bretagne Atlantique
(LMBA), UMR 6205, CNRS, Université de Bretagne Occidentale,
6 avenue Victor Le Gorgeu, 29238 Brest, France}
\email{luizlara@ime.unicamp.br}
\author[\'E. Loubeau]{\'Eric Loubeau}
\address{Laboratoire de Mathématiques de Bretagne Atlantique
(LMBA), UMR 6205, CNRS, Université de Bretagne Occidentale,
6 avenue Victor Le Gorgeu, 29238 Brest, France}
\email{eric.loubeau@univ-brest.fr}
\date{\today}
\hypersetup{
 pdftitle={On the Morse index of some Riemannian homogeneous fibrations},
 pdfauthor={Luiz Lara and \'Eric Loubeau},
 pdfsubject={Harmonic maps, Morse index and Jacobi spectra on homogeneous spaces},
 pdfkeywords={harmonic maps, Morse index, Jacobi operator, Hopf fibrations,
 Berger metrics, homogeneous spaces, Casimir operators,
 Aloff-Wallach spaces, Clifford hypersurfaces}
}

\begin{document}

\begin{abstract}
In this paper, we use the technique developed by the first named author in \cite{LaraComplexHopf} to compute the Morse-index of some Riemannian homogeneous fibrations.
For the quaternionic (resp. octonionic) Hopf fibration, we compute the
index and nullity, where the metric on the $\bbS^7$ (resp. $\bbS^{15}$) is the quaternionic (resp. octonionic) Hopf--Berger metric. In particular, for the round metric on the domain, the Morse-index is $45$ (resp. $153$). We also give a spectral formula for the Aloff--Wallach projections. In particular, we obtain that the Konishi $\R P^3$-fibration 
$N_{1,1}\to\C P^2$ is stable throughout the
two-parameter metric family, while the unit tangent sphere bundle $T^1(\C P^2)\to \C P^2$ is stable under some conditions on the metric. We also compute the index and nullity of
the Clifford hypersurface in the fiber-squashing range of the
quaternionic Berger family.
Finally, for the identity map on the quaternionic Hopf--Berger spheres, we determine its Morse--index by an explicit sequence of
transition parameters.
\end{abstract}
\maketitle

\tableofcontents

\section{Introduction}\label{sec:introduction}
\subsection{Context}
The Morse index of a harmonic map is the dimension of the
negative eigenspace of its Jacobi operator, and measures the number
of independent unstable directions for the Dirichlet energy.
For a nonconstant smooth harmonic map $\psi$ from the unit round
sphere $\bbS^m$, with $m\geq3$, El Soufi proved the sharp bound
\[
 \operatorname{Ind}(\psi)
 \geq\operatorname{Ind}(\operatorname{id}_{\bbS^m})=m+1;
\]
see~\cite[Theorem~1.1]{ElSoufi1995}. This leads to a natural question of chatacterizing the harmonic maps, if any, at which this minimal value for the index is attained.

For the round complex Hopf fibration
$\pi:\bbS^{2n+1}\to\C P^n$, this bound is attained: its index is
$2n+2$. For $n=1$, Urakawa~\cite{Urakawa1987-so} computed its index and
nullity as $4$ and $8$, respectively; this calculation was also
recovered by Loubeau and Oniciuc~\cite{LoubeauOniciuc2007}.
Rivi\`ere subsequently proved that every nonconstant harmonic map
$\bbS^3\to\bbS^2$ of index at most $4$ factors, up to an isometry
of the domain, through the Hopf fibration followed by a holomorphic
self-map of $\C P^1$~\cite{Riviere2023-co}. This result has been  generalized recently to harmonic maps $\bbS^{2n+1}\to\C P^n$ of index $2n+2$, see~\cite{LaraComplexHopf,GutweinLanglais2026,Surkov2026,ChenWangZhang2026}. They first show that the index of the complex Hopf map $\bbS^{2n+1}\to\C P^n$ is in fact $2n+2$ and that any other harmonic map with this index is, up to isometries of the sphere, the Hopf map composed with a holomorphic self-map of $\C P^n$. Also, in~\cite{LaraComplexHopf}, is presented a technique to compute the spectrum of the Jacobi operator of the Riemannian homogeneous fibration $(G/H,g)\to (G/K,h)$, given by $xH\mapsto xK$, where the domain is a geodesic-orbit space and the target is a symmetric space. 

These results suggest asking whether the other Hopf fibrations
also attain the universal lower bound and admit an analogous
equality-case characterisation. Our computations on this paper show a different behaviour. For the unit round domain metrics, the quaternionic
and octonionic Hopf projections
\[
 \bbS^7\longrightarrow\bbS^4,
 \qquad
 \bbS^{15}\longrightarrow\bbS^8
\]
have indices $45$ and $153$, respectively, rather than the lower
bounds $8$ and $16$ supplied by El Soufi's theorem.
Thus the direct analogue of the complex least-index equality case
does not describe these Hopf maps. This does not establish the
existence of harmonic maps between these same round spheres with
smaller indices. Determining the least Morse index among nonconstant
smooth harmonic maps $\bbS^7\to\bbS^4$ or
$\bbS^{15}\to\bbS^8$, and describing the maps which realise it,
remain open problems.

Changing the domain metric leads to a different question. The Hopf
projections remain harmonic when the fibres are squashed while the
horizontal metric is fixed. We study how their indices change under
this deformation, whether the resulting values attain known lower
bounds for maps into spheres, and whether the identity continues
to provide a lower bound for harmonic maps from the same domain.
The homogeneous method developed in~\cite{LaraComplexHopf} makes
these questions accessible: it expresses the Jacobi operator in
terms of Casimir operators and reduces its spectrum to branching
rules, tensor products and fixed subspaces. Here we apply it to
the two Hopf fibrations above and to the Aloff--Wallach projections,
and use related Casimir calculations for the Clifford hypersurface
and the identity on quaternionic Hopf--Berger spheres.

\subsection{Main results}
For either Hopf fibration, write
\[
 g_t=g_H+t\,g_V,\qquad 0<t\leq1,
\]
where $g_H$ and $g_V$ are the horizontal and vertical parts of the
corresponding unit round metric. Thus $t=1$ is the round metric and
$0<t<1$ is the fiber-squashing range. Denote the quaternionic and
octonionic projections, with domain metric $g_t$, by
$\phi_{\mathbb H,t}$ and $\phi_{\mathbb O,t}$, respectively;
the round target metrics are fixed so that the projections are
Riemannian submersions. We retain the inverse fibre-scale parameter
$\epsilon$ in the subsequent calculations, with
$t=(2\epsilon^2)^{-1}$ in the quaternionic case and
$t=\epsilon^{-2}$ in the octonionic case.

Theorems~\ref{thm:quaternionic-index} and
\ref{thm:octonionic-index} give
\[
 \operatorname{Ind}(\phi_{\mathbb H,t})=
 \begin{cases}
 45,&\frac34<t\leq1,\\
 13,&\frac38<t\leq\frac34,\\
 5,&0<t\leq\frac38,
 \end{cases}
 \qquad
 \operatorname{Ind}(\phi_{\mathbb O,t})=
 \begin{cases}
 153,&\frac7{16}<t\leq1,\\
 25,&\frac7{20}<t\leq\frac7{16},\\
 9,&0<t\leq\frac7{20}.
 \end{cases}
\]
We also determine their full Jacobi spectra, including
multiplicities, and their nullities at every parameter.
In particular, the quaternionic projection has indices $45$, $13$
and $5$ for the round, normal homogeneous and Jensen metrics,
respectively. The octonionic projection has index $153$ for the
round metric and index $9$ for both the normal homogeneous and
non-round Einstein metrics. The proofs identify all representation
types which can contribute nonpositive eigenvalues and show that
the remaining spectrum is strictly positive.

The quaternionic minimum $5$ has a direct variational meaning.
El Soufi~\cite[Theorem~2.5]{ElSoufi93} proved that a harmonic map
$\psi:(M,g)\to\bbS^4$ from a compact manifold has index at least
$5$ if its stress-energy tensor
\[
 S_\psi=\tfrac12|d\psi|^2g-\psi^*h
\]
is nonnegative everywhere and positive definite somewhere.
For the quaternionic Hopf Riemannian submersion, one has
$\tfrac12|d\phi_{\mathbb H,t}|^2=2$ and
\[
 S_{\phi_{\mathbb H,t}}=g_H+2t\,g_V>0.
\]
Consequently, for $0<t\leq3/8$, the Hopf projection attains this
lower bound. Squashing the fibres therefore produces maps of
index $5$ within this stress-energy class, including the projection
from Jensen's non-round Einstein sphere. The related rigidity
question for equality is discussed in
Remark~\ref{rem:quaternionic-index-rigidity}.

El Soufi's round-sphere comparison with the identity motivates the
calculation in the final section. Let
$\iota_t=\operatorname{id}_{(\bbS^7,g_t)}$ denote the identity,
with the same metric on its domain and target. Tanno studied
the Jacobi operator of the identity for quaternionic deformations
of round spheres~\cite{Tanno1993Identity}. In
Theorem~\ref{thm:identity-Berger-index}, we determine its Morse
index for the whole fiber-squashing range by an explicit sequence
of transition parameters and their jump multiplicities.
The behaviour is opposite to that of the Hopf projection:
$\operatorname{Ind}(\iota_t)$ tends to infinity as $t\to0$,
whereas $\operatorname{Ind}(\phi_{\mathbb H,t})$ eventually equals $5$.
In particular,
\[
 \operatorname{Ind}(\phi_{\mathbb H,t})=5
 <\operatorname{Ind}(\iota_t),\qquad 0<t\leq\frac38.
\]
At Jensen's parameter $t=1/5$, the two indices are $5$ and $13$.
Thus the inequality comparing every nonconstant harmonic map
with the identity of its domain, valid for round spheres, does
not extend to these homogeneous Berger metrics.

The same spectral method also applies to the Aloff--Wallach
projections
\[
 \phi_{k,l}:N_{k,l}=\SU(3)/H_{k,l}\longrightarrow\C P^2,
 \qquad
 H_{k,l}=\{\diag(z^k,z^l,z^{-k-l}):z\in\U(1)\},
\]
for coprime integers $k,l$, not both zero.
Theorem~\ref{thm:AW-index} gives their index and nullity for the
two-parameter homogeneous family $g_{\epsilon,\delta}$ by a finite
sum over eight possible intermediate representation types.
Two geometric specialisations are particularly useful.
The Konishi projection $\phi_{1,1}$, with fibre
$\mathrm{SO}(3)\cong\mathbb{R}P^3$, is stable with nullity $16$
throughout this family; see Corollary~\ref{cor:AW11-index}.
The natural unit tangent bundle projection
$T^1\C P^2\to\C P^2$, corresponding to $\phi_{1,-2}$ and having
fibre $\bbS^3$, is stable exactly when $\epsilon^2\geq\delta^2$;
see Corollary~\ref{cor:AW1m2-index}. This is the fiber-squashing
regime relative to the normal homogeneous metric. Outside it, its
index takes the values $2$, $18$ and $34$, with all transition
parameters and exceptional nullities given explicitly.

We also compute, using Casimir operators, the spectrum governing
normal variations of the Clifford hypersurface
$\bbS^3(1/\sqrt2)\times\bbS^3(1/\sqrt2)\subset\bbS^7$
throughout the quaternionic fiber-squashing range.
Here the Morse index is the index of the area functional under
normal variations. It equals $9$ for $1/3<t\leq1$ and $1$ for
$0<t\leq1/3$; in particular, the Clifford hypersurface in Jensen's
sphere has index $1$ and nullity $4$.
This calculation contributes to the study of variational indices
associated with non-round $\Sp(2)$-invariant metrics on $\bbS^7$.
While this manuscript was in preparation, Anna
Siffert~\cite{Siffert2026HopfSquashing} obtained these index and
nullity values as part of a broader study of Clifford hypersurfaces
under Hopf squashing. Our treatment provides a Casimir-operator
derivation of this special case.

\subsection{Organisation}
Section~\ref{sec:homogeneous-spectral-framework} fixes the notation
and recalls the homogeneous Jacobi spectral formula from
\cite[Theorem~2.2]{LaraComplexHopf}, to which we refer for the proofs.
Section~\ref{sec:Hopf-fibrations} applies it to the quaternionic and
octonionic Hopf fibrations and determines their index and nullity
throughout the fiber-squashing range.
Section~\ref{sec:aloff-wallach} treats the Aloff--Wallach projections,
including the Konishi and unit tangent bundle specialisations.
Section~\ref{sec:Clifford-hypersurface} gives the Casimir calculation
for the Clifford hypersurface. Finally,
Section~\ref{sec:Berger-identity} computes the index of the identity
on quaternionic Hopf--Berger spheres and establishes its divergence
under fibre collapse.

\section{Jacobi operator on a homogeneous fibration}\label{sec:homogeneous-spectral-framework}
We recall the geometric setting and representation-theoretic notation
needed to apply~\cite[Theorem~2.2]{LaraComplexHopf}. The main result of
this section, Theorem~\ref{thm:J-eigenvalues}, expresses the Jacobi
spectrum of a homogeneous Riemannian fibration in terms of Casimir
eigenvalues, branching rules and $H$-fixed subspaces. All results in
this section are recalled without proof; their proofs are given in
\cite[Section~2]{LaraComplexHopf}.

\subsection{Riemannian homogeneous fibrations}\label{ssec:RHF}
Let $G$ be a compact connected Lie group and let $H\subset K\subset G$
be closed connected subgroups. We denote their Lie algebras by $\fh$, $\fk$ and
$\fg$, and put $s=\dim G$, $m=\dim(G/H)$ and $n=\dim(G/K)$.
Fix an $\Ad(G)$-invariant inner product $\bangle{\cdot,\cdot}$ on $\fg$.
Let $\fp$ be the orthogonal complement of $\fh$ in $\fk$, and let
$\fq$ be the orthogonal complement of $\fk$ in $\fg$. Thus
\begin{equation}\label{eq:cartan-decomp}
 \fg=\fk\oplus\fq=\fh\oplus\fp\oplus\fq,
 \qquad \fm=\fp\oplus\fq.
\end{equation}
The subspaces $\fp$ and $\fq$ are the vertical and horizontal
subspaces of $G/H\to G/K$ at $eH$, respectively, and we identify
$\fm$ with $T_{eH}(G/H)$. The corresponding $G$-invariant
distributions are denoted by $\mathcal V$ and $\mathcal H$.
For $X\in\fm$, the subscripts in $X_\fp$ and $X_\fq$ denote
the orthogonal projections onto these summands. In particular,
\begin{equation}\label{eq:invariance-cartan}
 \Ad(K)\fq=\fq,\qquad \Ad(H)\fp=\fp,\qquad \Ad(H)\fm=\fm.
\end{equation}

Given $\epsilon,\delta>0$, consider the inner products
\begin{equation}\label{eq:homogeneous-metric}
 \bangle{X,Y}_{\epsilon,\delta}
  =\epsilon^{-2}\bangle{X_\fp,Y_\fp}
     +\delta^{-2}\bangle{X_\fq,Y_\fq},
 \qquad X,Y\in\fm,
\end{equation}
and $\delta^{-2}\bangle{\cdot,\cdot}$ on $\fq$.
We denote the induced $G$-invariant metrics on $G/H$ and $G/K$ by
$g_{\epsilon,\delta}$ and $h_\delta$, respectively. The projection
\[
 \phi:(G/H,g_{\epsilon,\delta})\longrightarrow(G/K,h_\delta),
 \qquad \phi(xH)=xK,
\]
is a Riemannian submersion with totally geodesic fibres, and hence is
harmonic.

We write $E=\phi^{-1}(T(G/K))$ and denote its pullback connection by
$\nabla^\phi$. With the curvature convention
$R^h(X,Y)Z=\nabla^h_X\nabla^h_YZ-\nabla^h_Y\nabla^h_XZ
-\nabla^h_{[X,Y]}Z$, the Jacobi operator is
\begin{equation}\label{eq:Jacobi-definition}
 \cJ_\phi V=(\nabla^\phi)^\ast\nabla^\phi V
    -\sum_{i=1}^m R^h(V,d\phi(e_i))d\phi(e_i),
\end{equation}
where $\{e_i\}$ is a local $g_{\epsilon,\delta}$-orthonormal frame.
It is a self-adjoint elliptic operator. We use
$\mathrm{Ind}(\phi)$ for the number of its negative eigenvalues,
counted with real multiplicity, and
$\mathrm{Null}(\phi)=\dim_\R\ker\cJ_\phi$ for its nullity.
Scalar Laplacians have the nonnegative convention
$\Delta_g=-\operatorname{div}_g\nabla^g$.

The bundle $E$ is naturally identified with $G\times_H\fq$. Consequently,
its smooth sections correspond to the equivariant functions
\begin{equation}\label{def:H-eq-functions}
 C^\infty_H(G,\fq)=
 \{f\in C^\infty(G,\fq):f(xh)=\Ad(h^{-1})f(x)
                  \text{ for }x\in G,\ h\in H\}.
\end{equation}
For $x\in G$, let $\tau_x$ denote left translation on the homogeneous
spaces, and let $\Gamma(E)$ denote the space of smooth sections of $E$.
The identification is the $G$-equivariant isomorphism
\begin{equation}\label{eq:Psi-section-identification}
 \Psi:C^\infty_H(G,\fq)\longrightarrow\Gamma(E),
 \qquad \Psi(f)(xH)=(\tau_x)_\ast f(x).
\end{equation}
If a local section of $G\to G/H$ is understood, we also write
\begin{equation}\label{eq:good-frame}
 X^\ast_{xH}=(\tau_x)_\ast X,\qquad X\in\fm.
\end{equation}

We write $\Ad$ for the adjoint group action and
$\ad(X)Y=[X,Y]$ for its differential.
For $a\in G$, let $R_a(x)=xa$ denote right translation.
For $X\in\fg$, let
$Xf=\left.\frac{d}{dt}\right|_{t=0}f\circ R_{\exp(tX)}$.
For $X\in\fk$, put $\Theta(X)f=Xf+\ad(X)f$.
These first-order operations need not individually preserve
$H$-equivariance, but the invariant quadratic sums used below do.
Choose a $\bangle{\cdot,\cdot}$-orthonormal basis $\{Y_i\}_{i=1}^s$
adapted to $\fg=\fq\oplus\fp\oplus\fh$. Thus $Y_1,\ldots,Y_n$ lie
in $\fq$, $Y_{n+1},\ldots,Y_m$ in $\fp$, and the remaining vectors in
$\fh$. Define
\begin{equation}\label{eq:jacobian-pieces}
 \begin{aligned}
 D_G&=\sum_{i=1}^sY_i^2,&
 D_K&=\sum_{i=n+1}^sY_i^2,\\
 D_G^\ad&=\sum_{i=1}^s\ad(Y_i)^2,&
 D_K^\ad&=\sum_{i=n+1}^s\ad(Y_i)^2,&
 D_\Theta&=\sum_{i=n+1}^s\Theta(Y_i)^2.
 \end{aligned}
\end{equation}
The adjoint sums are first taken on $\fg$; their difference preserves
$\fq$.

\begin{theorem}[{\cite[Theorem~2.1]{LaraComplexHopf}}]\label{thm:Urakawa}
Assume that $(G/H,g_{\epsilon,\delta})$ is a geodesic orbit space
and that $G/K$ is a symmetric space. Under the
identification~\eqref{eq:Psi-section-identification}, the Jacobi
operator of $\phi$ is
\begin{equation}\label{eq:Jacobi-Casimir-operator}
 \cJ_\phi
   =\delta^2(D_K-D_G+D_G^\ad-D_K^\ad)-\epsilon^2D_\Theta.
\end{equation}
\end{theorem}
\subsection{Notations for representation theory}
We use $\N_0=\{0,1,2,\ldots\}$ and write $\mathbf1_A$ for the
indicator of a condition or set $A$.
Fix maximal tori $T_K\subset T_G$ and compatible positive root
systems for $G$ and $K$ in the duals of their torus Lie algebras.
We write
$\Lambda_G^+$ and $\Lambda_K^+$ for the sets of dominant integral weights,
and $\rho_G$ and $\rho_K$ for the half-sums of their positive roots.
The inner product on weights is induced by the chosen invariant
inner product, restricted to the corresponding torus Lie algebra.
For a semisimple group $L$, the notation $\varpi_i^L$ denotes
the fundamental weights: if $\alpha_j^L$ are the chosen simple
roots, then
$2\bangle{\varpi_i^L,\alpha_j^L}/\bangle{\alpha_j^L,\alpha_j^L}
=\delta_{ij}$, where $\delta_{ij}$ is the Kronecker symbol.
When $K$ has a central torus, the coordinates of its full weights,
including the central characters, will be specified in the example.
The superscript on a weight or representation identifies its group;
the coordinate bases and index ranges are fixed separately in each example.
For $\lambda\in\Lambda_G^+$, let
$\pi_\lambda^G:G\to\GL(\cH_\lambda^G)$ denote the irreducible complex
representation with highest weight $\lambda$; we use the analogous
notation for $K$.
The superscript $\ast$ on a representation denotes its contragredient,
and on a vector space its complex linear dual; on differential
operators, as in $(\nabla^\phi)^\ast$, it denotes the formal adjoint.
For a real vector space $W$, put $W_\C=W\otimes_\R\C$.
If $U,V$ are complex $H$-modules, $\Hom_H(U,V)$ denotes the space
of complex linear intertwiners, and
$U^H=\{u\in U:h\cdot u=u\text{ for all }h\in H\}$.
The tensor product of representations of the same group carries
the diagonal action. For representations of two different group
factors, $\boxtimes$ denotes the external tensor product, with
each factor acting on its own vector-space factor.
Let $\pi:G\to\GL(V)$ be a finite-dimensional complex
representation, and let $\pi_\ast:\fg\to\gl(V)$ be its
differential; we also write $d\pi=\pi_\ast$.
Given a $\bangle{\cdot,\cdot}$-orthonormal basis
$\{Y_i\}_{i=1}^s$ of $\fg$, the Casimir operator of $\pi$ is the
endomorphism of $V$ defined by
\[
 \Cas_G^\pi=\sum_{i=1}^s\pi_\ast(Y_i)^2.
\]
This operator is independent of the chosen orthonormal basis and
commutes with $\pi(g)$ for every $g\in G$. For the irreducible representation
$\pi_\lambda^G$,
\[
 \Cas_G^{\pi_\lambda^G}=-c_G(\lambda)\id_{\cH_\lambda^G},
\]
where
\[
 c_G(\lambda)\coloneqq\bangle{\lambda+2\rho_G,\lambda},
\]
and the inner product on weights is induced by
$\bangle{\cdot,\cdot}$ on $\fg$. We call $c_G(\lambda)$ the
Casimir constant; the eigenvalue of $\Cas_G^{\pi_\lambda^G}$ is
$-c_G(\lambda)$.
The operators $\Cas_K^\pi$ and the constants $c_K(\mu)$ are defined
analogously, using the restriction of the same inner product to $\fk$.

We write $\Ad_\fq(k)=\Ad(k)\rvert_\fq$ for the real
$K$-isotropy representation. Recall that a real irreducible representation may be of real,
complex or quaternionic type. We denote by
$(\pi_\fq,\cV_\fq)$ the associated complex irreducible
representation when $\Ad_\fq$ is irreducible, as in the
irreducible symmetric-space setting below. The type always
refers to the $K$-representation; restriction to $H$ may change it.
Thus
\[
 \fq_\C\cong
 \begin{cases}
  \cV_\fq,&\text{if $\Ad_\fq$ is of real type},\\
  \cV_\fq\oplus\cV_\fq^\ast,
      &\text{if $\Ad_\fq$ is of complex type},\\
  \cV_\fq\oplus\cV_\fq,
      &\text{if $\Ad_\fq$ is of quaternionic type}.
 \end{cases}
\]
Let $L^2_H(G,\fq_\C)$ be the Hilbert completion of
$C^\infty_H(G,\fq_\C)$, using Haar measure and the invariant
Hermitian inner product on $\fq_\C$.
Its Peter--Weyl decomposition is
\begin{equation}\label{eq:Peter-Weyl-pullback}
 L^2_H(G,\fq_\C)
  \cong\widehat\bigoplus_{\lambda\in\Lambda_G^+}
       \cH_\lambda^G\otimes\Hom_H(\cH_\lambda^G,\fq_\C).
\end{equation}
Here $\widehat\bigoplus$ is the Hilbert direct sum; smooth sections
are the smooth vectors for the $G$-action in this completion.
Assume from now on that $G$ is simple and $G/K$ is an irreducible
symmetric space.
Let $\mu_\fg$ be the highest weight of the complex adjoint
representation, and let $\mu_\fq$ be the highest weight of the
chosen representation $\pi_\fq$. The scalar
$c_K(\mu_\fq)-c_G(\mu_\fg)$ is the isotropy contribution to the
spectral formula. Following~\cite{LaraComplexHopf}, we use the
following notation.
\begin{definition}\label{def:C-lambda}
For $\lambda\in\Lambda_G^+$ and $\mu_1,\mu_2\in\Lambda_K^+$, set
\[
 C_\lambda(\mu_1,\mu_2)
  \coloneqq\delta^2\Big(c_G(\lambda)-c_K(\mu_1)
       +c_K(\mu_\fq)-c_G(\mu_\fg)\Big)
       +\epsilon^2c_K(\mu_2).
\]
The dependence on $\epsilon,\delta$ is suppressed in this notation.
\end{definition}
\begin{theorem}[{\cite[Theorem~2.2]{LaraComplexHopf}}]\label{thm:J-eigenvalues}
Let $\phi:(G/H,g_{\epsilon,\delta})\to(G/K,h_\delta)$ be the
homogeneous fibration above, with $G$ simple and $G/K$ irreducible
symmetric. Decompose
\[
 \pi_\lambda^G\rvert_K=\bigoplus_i\pi_{\lambda_i}^K,
\]
and, for every $i$, decompose
\[
 (\pi_{\lambda_i}^K)^\ast\otimes\pi_\fq
    =\bigoplus_j\pi_{\lambda_{ij}}^K,
 \qquad
 (\pi_{\lambda_i}^K)^\ast\otimes\pi_\fq^\ast
    =\bigoplus_j\pi_{\bar\lambda_{ij}}^K,
\]
listing repeated summands separately. The two lists of final
weights are indexed independently: the bar in $\bar\lambda_{ij}$
labels the second tensor product, not conjugation of an individual weight.
Put
\[
 d_\lambda=\dim_\C\cH_\lambda^G,
 \qquad
 m_\nu=\dim_\C(\cH_\nu^K)^H.
\]
The eigenvalues and their real multiplicities are as follows:
\begin{center}
\begin{tabular}{lcc}
\toprule
Type of $\Ad_\fq$ & Eigenvalues & Multiplicities\\
\midrule
Real
 & $C_\lambda(\lambda_i,\lambda_{ij})$
 & $m_{\lambda_{ij}}d_\lambda$\\
Complex
 & $C_\lambda(\lambda_i,\lambda_{ij})$,
   $C_\lambda(\lambda_i,\bar\lambda_{ij})$
 & $m_{\lambda_{ij}}d_\lambda$,
   $m_{\bar\lambda_{ij}}d_\lambda$\\
Quaternionic
 & $C_\lambda(\lambda_i,\lambda_{ij})$
 & $2m_{\lambda_{ij}}d_\lambda$\\
\bottomrule
\end{tabular}
\end{center}
\end{theorem}
\section{Jacobi spectra of Hopf fibrations}\label{sec:Hopf-fibrations}
We now apply Theorem~\ref{thm:J-eigenvalues} to the Hopf fibrations.
For each fibration, let $g_{\mathrm{rd}}$ be the unit round metric
and write $g_H=g_{\mathrm{rd}}|_{\mathcal H}$ and
$g_V=g_{\mathrm{rd}}|_{\mathcal V}$ for its horizontal and vertical
parts, extended by zero on the complementary distributions.
To specify the fiber-squashing range, we compare our
metrics with the convention
\begin{equation}\label{eq:Hopf-short-fibre-convention}
 g_t=g_H+t\,g_V,\qquad 0<t\leq1.
\end{equation}
Thus $t$ is the squared fibre-length scale, $t=1$ is the round
metric, and $0<t<1$ squashes the fibers. The parameter $\epsilon$
in the Casimir formula remains our parameter in the spectral
calculations; its relation to $t$ is specified in each example.
We first recall the complex case from~\cite{LaraComplexHopf} to fix
the comparison of metric parameters and real multiplicities, and
then carry out the quaternionic and octonionic calculations.

\subsection{Comparison with the complex Hopf fibration}
For comparison, we recall the corresponding calculation for the
complex Hopf map, following the normalisation and results
of~\cite{LaraComplexHopf}.
Let
\[
 G=\SU(n+1),\qquad H=\SU(n),\qquad
 K=\mathrm{S}(\U(1)\times\U(n)),
\]
where $H$ acts on the last $n$ coordinates.
Thus $G/H\cong\bbS^{2n+1}$ and $G/K\cong\C P^n$, and
$\phi$ is the quotient projection. We choose
$\bangle{X,Y}=-\frac12\tr(XY)$ on $\fg$.
We fix $\delta=1$ and write $g_\epsilon=g_{\epsilon,1}$ and $h=h_1$.
If $g_{\mathrm{rd}}$ is the unit round metric, $\xi$ its unit Hopf
vector field, and $\eta=g_{\mathrm{rd}}(\xi,\cdot)$, then
\begin{equation}\label{eq:complex-round-normalisation}
 g_\epsilon=g_{\mathrm{rd}}
      +\left(\frac{n+1}{2n\epsilon^2}-1\right)\eta\otimes\eta.
\end{equation}
In the convention~\eqref{eq:Hopf-short-fibre-convention},
\[
 t=\frac{n+1}{2n\epsilon^2},\qquad
 \epsilon^2\geq\frac{n+1}{2n}.
\]
The round parameter is $\epsilon^2=(n+1)/(2n)$.
Here $\Ad_\fq$ is of complex type, so both conjugate summands in
$\fq_\C$ must be included when counting real multiplicities.

For $n\geq3$, the only eigenvalue which can be negative in this
range is
\[
 \frac{2n}{n+1}\epsilon^2-2n,
\]
with real multiplicity $2(n+1)$. The block with eigenvalue
$\frac{8n}{n+1}\epsilon^2-4$ has multiplicity $n(n+1)$;
it is positive for $\epsilon^2>\frac{n+1}{2n}$ and vanishes
at the round metric.
Consequently, the index and nullity are
\begin{align}
 \mathrm{Ind}(\phi)&=
 \begin{cases}
 2(n+1),&\frac{n+1}{2n}\leq\epsilon^2<n+1,\\
 0,&\epsilon^2\geq n+1,
 \end{cases}\label{eq:complex-index-recalled}\\
 \mathrm{Null}(\phi)&=
 \begin{cases}
 n(3n+5),&\epsilon^2=\frac{n+1}{2n},\\
 2(n^2+3n+1),&\epsilon^2=n+1,\\
 2n(n+2),&\text{otherwise}.
 \end{cases}\label{eq:complex-nullity-recalled}
\end{align}
The branching reduction and the treatment of the two conjugate
isotropy summands are given in~\cite{LaraComplexHopf}.
This comparison highlights the distinction between the complex
case and the two real-type isotropy representations considered
here; the metric parameter and multiplicity conventions are the
same throughout.

\subsection{The quaternionic Hopf fibration}\label{subsec:quaternionic-Hopf}
Let $G=\Sp(2)$, $K=\Sp(1)\times\Sp(1)$ and
\[
 H=\{\diag(h,1):h\in\Sp(1)\}\subset K.
\]
The Hopf map is the homogeneous Riemannian submersion
$\phi:G/H\cong\bbS^7\to G/K\cong\H P^1\cong\bbS^4$.
We choose the $\Ad(G)$-invariant inner product
\[
 \bangle{X,Y}=-\frac12\Re\tr_{\H}(XY)
\]
on $\fg$. Here the trace is taken in the defining quaternionic
representation. We fix $\delta=1$ and write $g_\epsilon=g_{\epsilon,1}$.
With the notation of~\eqref{eq:cartan-decomp}, these metrics satisfy
\begin{equation}\label{eq:quaternionic-metric}
 g_\epsilon=\bangle{\cdot,\cdot}\rvert_{\fq}
       +\epsilon^{-2}\bangle{\cdot,\cdot}\rvert_{\fp},
 \qquad
 g_{\mathrm{rd}}=\bangle{\cdot,\cdot}\rvert_{\fq}
       +2\bangle{\cdot,\cdot}\rvert_{\fp}.
\end{equation}
Thus $g_{\mathrm{rd}}$ is the round metric of radius $1$ when
$\epsilon^2=\frac12$, and $h_1$ is the round metric of radius $\frac12$
on $\bbS^4$; see~\cite[Section 3.2]{Bettiol2022-vw}.
In the convention~\eqref{eq:Hopf-short-fibre-convention},
the parameters and their admissible ranges are
\begin{equation}\label{eq:quaternionic-t-epsilon}
 t=\frac{1}{2\epsilon^2},\qquad
 0<t\leq1\quad\Longleftrightarrow\quad
 \epsilon\geq\frac1{\sqrt2}.
\end{equation}
We henceforth restrict this example to that range.
The fibres are totally geodesic, so $\phi$ is harmonic throughout it.

We first fix the weight conventions. Let
$\mathbf i,\mathbf j,\mathbf k$ denote the standard imaginary
quaternionic units. A maximal torus of $G$ is
\[
 T=T_G=T_K=\{\diag(e^{\mathbf i\theta_1},e^{\mathbf i\theta_2}):
                   \theta_1,\theta_2\in\R\}.
\]
In the associated coordinates $e_1,e_2$, the positive roots of $G$ are
$e_1-e_2,e_1+e_2,2e_1,2e_2$, and the positive roots of $K$ are
$2e_1,2e_2$. In the fundamental-weight notation fixed in
Section~\ref{sec:homogeneous-spectral-framework},
we have
\[
 \varpi_1^G=e_1,\qquad \varpi_2^G=e_1+e_2,\qquad
 \varpi_1^K=e_1,\qquad \varpi_2^K=e_2,
 \qquad \rho_G=(2,1),\qquad \rho_K=(1,1).
\]
We write $\pi^G_{(a,b)}:G\to\GL(\cH^G_{(a,b)})$ for the
irreducible complex representation of highest weight
$(a,b)=ae_1+be_2$, where $a\geq b\geq0$ are integers.
Equivalently, this weight is
$(a-b)\varpi_1^G+b\varpi_2^G$: the pair $(a,b)$ gives
coordinate coefficients, rather than Dynkin labels.
Similarly, $\pi^K_{(r,s)}:K\to\GL(\cH^K_{(r,s)})$ denotes the
irreducible representation of highest weight $(r,s)$,
where $r,s\geq0$ are integers.
We abbreviate $c_G((a,b))$ and $c_K((r,s))$ to
$c_G(a,b)$ and $c_K(r,s)$, respectively.
If $\pi_r^{\Sp(1)}$ denotes the irreducible complex
$\Sp(1)$-representation of highest weight $r$, on
$\cH_r^{\Sp(1)}$ of dimension $r+1$, then
$\pi^K_{(r,s)}=\pi_r^{\Sp(1)}\boxtimes\pi_s^{\Sp(1)}$.
Here the external tensor product acts explicitly by
\[
 (\pi_r^{\Sp(1)}\boxtimes\pi_s^{\Sp(1)})(k_1,k_2)(u\otimes v)
 =\pi_r^{\Sp(1)}(k_1)u\otimes\pi_s^{\Sp(1)}(k_2)v.
\]
The induced inner product on weights is
$\bangle{(a,b),(a',b')}=2aa'+2bb'$. Consequently,
\begin{align}
 c_G(a,b)&=\bangle{(a,b)+2\rho_G,(a,b)}
          =2\big(a(a+4)+b(b+2)\big),\label{eq:cas-Sp2}\\
 c_K(r,s)&=\bangle{(r,s)+2\rho_K,(r,s)}
          =2\big(r(r+2)+s(s+2)\big).\label{eq:cas-Sp1xSp1}
\end{align}
The highest weights of the complexified adjoint and isotropy representations are
$\mu_\fg=(2,0)$ and $\mu_\fq=(1,1)$, respectively.
The isotropy representation $\Ad_\fq$ is of real type, and its
associated complex representation is $\pi_\fq=\pi^K_{(1,1)}$.
We now describe the final highest weights $\lambda_{ij}$ of
Theorem~\ref{thm:J-eigenvalues}.

\begin{lemma}
The isotropy contribution to $C_\lambda$ is
\begin{equation}\label{eq:quaternionic-Cp}
 c_K(\mu_\fq)-c_G(\mu_\fg)
   =c_K(1,1)-c_G(2,0)=12-24=-12.
\end{equation}
From the Weyl dimension formula, we obtain
\begin{equation}\label{eq:quaternionic-dimension}
 d_{a,b}\coloneqq\dim_\C \cH^G_{(a,b)}
  =\frac{(a-b+1)(a+b+3)(a+2)(b+1)}{6}.
\end{equation}
\end{lemma}
We now determine the $K$-components contributing to
$\Hom_H(\cH^G_{(a,b)},\fq_\C)$, starting with the branching rule
for $G\downarrow K$.

\begin{theorem}[Branching rule, {\cite[Chapter~8]{GoodmanWallach2009}}]
\label{thm:quaternionic-branching}
Let $G=\Sp(2)$ and $K=\Sp(1)\times\Sp(1)$, embedded as the
block-diagonal subgroup. For integers $a\geq b\geq0$, the
restriction of $\pi^G_{(a,b)}$ to $K$ decomposes as
\begin{equation}\label{eq:quaternionic-branching}
 \pi^G_{(a,b)}\rvert_K
 \cong\bigoplus_{\substack{r,s\in\mathbb Z_{\geq0}\\
 a+b\geq r+s\geq a-b\geq|r-s|\\
 a+b-r-s\in2\mathbb Z}}\pi^K_{(r,s)}.
\end{equation}
Every summand occurs with multiplicity one.
\end{theorem}
This is the $\Spin(5)\downarrow\Spin(4)$ interlacing rule, under
$\Sp(2)\cong\Spin(5)$ and
$\Sp(1)\times\Sp(1)\cong\Spin(4)$.
\begin{lemma}\label{lem:quaternionic-tensor-product}
The weights $\lambda\in\Lambda_G^+$ satisfying
\[
 \Hom_H(\cH^G_\lambda,\fq_\C)\neq0
\]
are exactly
\begin{equation}\label{eq:quaternionic-contributing-weights}
 \lambda=ae_1+be_2,\qquad a\geq b\geq0,\quad a\geq1.
\end{equation}
For such a weight, put $x=a-b$. The $K$-irreducible components
whose restrictions to $H$ can contribute to this space are
\begin{equation}\label{eq:quaternionic-relevant-K-components}
 \pi^G_{(a,b)}\rvert_K
 \cong\tau_{a,b}
 \oplus\begin{cases}
  \pi^K_{(1,x-1)},&x\geq1,\\
  0,&x=0,
 \end{cases}
 \oplus\begin{cases}
  \pi^K_{(1,x+1)},&b\geq1,\\
  0,&b=0.
 \end{cases}
\end{equation}
Here $\tau_{a,b}$ is a $K$-representation on $\cH_{\tau_{a,b}}$ satisfying
$\Hom_H(\cH_{\tau_{a,b}},\fq_\C)=0$, and every displayed
irreducible summand has multiplicity one.

In the notation of Theorem~\ref{thm:J-eigenvalues}, choose
$\pi_\fq=\pi^K_{(1,1)}$. For integers $r,s\geq0$, the tensor
products decompose into $K$-irreducible representations as follows:
\begin{equation}\label{eq:quaternionic-tensor-product}
 \begin{aligned}
 (\pi^K_{(r,s)})^\ast\otimes\pi_\fq
 &\cong\pi^K_{(r+1,s+1)}\oplus\pi^K_{(r+1,s-1)}\\
 &\qquad\oplus\pi^K_{(r-1,s+1)}\oplus\pi^K_{(r-1,s-1)}.
 \end{aligned}
\end{equation}
Every summand has multiplicity one, and summands with negative
indices are omitted. The $H$-fixed part is nonzero if and only
if $r=1$. For the relevant intermediate type $\pi^K_{(1,s)}$,
the only final types with nonzero $H$-fixed subspaces are
$\pi^K_{(0,s+1)}$ and, when $s\geq1$, $\pi^K_{(0,s-1)}$, with
\begin{equation}\label{eq:quaternionic-final-types}
 \begin{aligned}
  (\cH^K_{(0,s+1)})^H&\cong\C^{s+2},&&s\geq0,\\
  (\cH^K_{(0,s-1)})^H&\cong\C^s,&&s\geq1.
 \end{aligned}
\end{equation}
The $H$-fixed parts of all the other summands vanish.
\end{lemma}
\begin{proof}
The restriction of $\fq_\C$ to $H=\Sp(1)\times\{1\}$ is the
sum of two copies of $\pi_1^{\Sp(1)}$, while
$\pi^K_{(r,s)}\rvert_H$ is the sum of $s+1$ copies of
$\pi_r^{\Sp(1)}$. Schur's lemma therefore gives
\[
 \Hom_H(\cH^K_{(r,s)},\fq_\C)\neq0
 \quad\Longleftrightarrow\quad r=1.
\]
We apply Theorem~\ref{thm:quaternionic-branching} with first
$K$-weight equal to $1$. Writing the second weight as $s$,
the interlacing and parity conditions give
\[
 |1-s|\leq x\leq1+s,\qquad 1+s\leq a+b,
 \qquad a+b-1-s\in2\mathbb Z.
\]
The first condition leaves $s=x-1,x,x+1$, and parity excludes
$s=x$. The choice $s=x-1$ is admissible exactly when $x\geq1$;
the choice $s=x+1$ is admissible exactly when $b\geq1$.
This proves~\eqref{eq:quaternionic-relevant-K-components}.
At least one of these components occurs for every $(a,b)\neq(0,0)$,
whereas neither occurs for $(a,b)=(0,0)$, proving
\eqref{eq:quaternionic-contributing-weights}.

Using self-duality and the Clebsch--Gordan formula for
$\Sp(1)\cong\SU(2)$~\cite[p.~339]{GoodmanWallach2009}, we obtain
\[
 \begin{aligned}
 (\pi^K_{(r,s)})^\ast\otimes\pi_\fq
 &\cong
 (\pi_r^{\Sp(1)}\boxtimes\pi_s^{\Sp(1)})^\ast
 \otimes(\pi_1^{\Sp(1)}\boxtimes\pi_1^{\Sp(1)})\\
 &\cong
 \big((\pi_r^{\Sp(1)})^\ast\otimes\pi_1^{\Sp(1)}\big)
 \boxtimes
 \big((\pi_s^{\Sp(1)})^\ast\otimes\pi_1^{\Sp(1)}\big)\\
 &\cong
 (\pi_{r+1}^{\Sp(1)}\oplus\pi_{r-1}^{\Sp(1)})
 \boxtimes
 (\pi_{s+1}^{\Sp(1)}\oplus\pi_{s-1}^{\Sp(1)}),
 \end{aligned}
\]
where $\pi_{-1}^{\Sp(1)}$ is understood to be zero.
Distributing the external tensor product over the direct sums
gives~\eqref{eq:quaternionic-tensor-product}, with each summand
occurring once.

Since $H$ acts only on the first factor of
$\cH^K_{(u,v)}=\cH_u^{\Sp(1)}\otimes\cH_v^{\Sp(1)}$, this
space has $H$-fixed vectors precisely when $u=0$.
In that case, the whole space is fixed and has dimension $v+1$.
In~\eqref{eq:quaternionic-tensor-product}, the first highest
weight is $r+1$ or $r-1$, so it can be zero exactly when $r=1$.
The surviving summands are then $(0,s+1)$ and $(0,s-1)$,
with the latter omitted when $s=0$. Their dimensions give
\eqref{eq:quaternionic-final-types}.
\end{proof}

Relabelling the second weight, we write the relevant intermediate
type as $\pi^K_{(1,r)}$. Applying Theorem~\ref{thm:J-eigenvalues}
and Lemma~\ref{lem:quaternionic-tensor-product}, we obtain the
following description of the entire spectrum.
\begin{proposition}\label{prop:quaternionic-spectrum}
Assume $\epsilon^2\geq\frac12$. For each $a\geq b\geq0$ and
each $r\geq0$ for which
$\pi^K_{(1,r)}\subset\pi^G_{(a,b)}\rvert_K$, the corresponding
eigenvalues of $\cJ_\phi$ are
\begin{equation}\label{eq:quaternionic-spectrum}
 \begin{aligned}
 C^\pm_{a,b;r}(\epsilon)
 &\coloneqq C_{(a,b)}\big((1,r),(0,r\pm1)\big)\\
 &=2\big(a(a+4)+b(b+2)-r(r+2)-9\big)
       +2\epsilon^2(r\pm1)(r\pm1+2).
 \end{aligned}
\end{equation}
The plus sign has multiplicity $d_{a,b}(r+2)$, and the minus sign,
which is present only for $r\geq1$, has multiplicity $d_{a,b}r$.
Multiplicities are added whenever different blocks have the same
numerical eigenvalue.

The only blocks which can be negative are
\begin{equation}\label{eq:quaternionic-negative-blocks}
 \begin{aligned}
 C^+_{1,0;0}(\epsilon)&=6\epsilon^2-8<0
       &&\Longleftrightarrow\quad\epsilon^2<\frac43,\\
 C^-_{1,1;1}(\epsilon)&=-8<0
       &&\text{for every }\epsilon^2\geq\frac12,\\
 C^-_{2,1;2}(\epsilon)&=6\epsilon^2-4<0
       &&\Longleftrightarrow\quad\epsilon^2<\frac23.
 \end{aligned}
\end{equation}
All the other blocks are nonnegative. Moreover,
\begin{equation}\label{eq:quaternionic-permanent-zero}
 C^-_{2,0;1}(\epsilon)=0
 \qquad\text{for every }\epsilon^2\geq\frac12,
\end{equation}
with real multiplicity $10$. Apart from this permanent zero block,
zero occurs only when
\begin{equation}\label{eq:quaternionic-zero-blocks}
 \begin{aligned}
 C^+_{1,0;0}(\epsilon)=0
       &\quad\Longleftrightarrow\quad\epsilon^2=\frac43,\\
 C^-_{2,1;2}(\epsilon)=0
       &\quad\Longleftrightarrow\quad\epsilon^2=\frac23,\\
 C^+_{1,1;1}(\epsilon)=0
       &\quad\Longleftrightarrow\quad\epsilon^2=\frac12.
 \end{aligned}
\end{equation}
\end{proposition}
\begin{proof}
For the intermediate type $\pi^K_{(1,r)}$,
$c_K(1,r)=6+2r(r+2)$. Substituting this and
\eqref{eq:quaternionic-Cp} in Theorem~\ref{thm:J-eigenvalues}, and
using the final types~\eqref{eq:quaternionic-final-types}, gives
\eqref{eq:quaternionic-spectrum}. The multiplicities follow from
the dimensions of their $H$-fixed subspaces
and~\eqref{eq:quaternionic-dimension}.
All computations are performed in the complexification of the
real bundle $\phi^{-1}(T\bbS^4)$. Since $\cJ_\phi$ is a real operator,
the complex dimension of each complexified eigenspace equals the
real dimension of its real eigenspace. Since $\Ad_\fq$ is of real
type, no additional factor of $2$ is introduced; see
\cite[Section~2.3]{LaraComplexHopf}.

For the sign analysis, put $x=a-b$. By
Lemma~\ref{lem:quaternionic-tensor-product}, the relevant
intermediate types have either $r=x+1$, with $b\geq1$, or
$r=x-1$, with $x\geq1$. Writing
$B=a(a+4)+b(b+2)-r(r+2)-9$, we obtain
\begin{equation}\label{eq:quaternionic-finite-reduction}
 \begin{aligned}
 r=x+1:\quad &B=2b(b+x+3)-12,\\
 r=x-1:\quad &B=2b(b+x+3)+4x-8.
 \end{aligned}
\end{equation}
Every block with $B>0$ is strictly positive. The first expression
is negative only for $(b,x)=(1,0),(1,1)$ and vanishes only for
$(b,x)=(1,2)$. The second expression is negative only for
$(b,x)=(0,1)$ and vanishes only for $(b,x)=(0,2)$.
The negative constant terms therefore give exactly the five blocks
\[
 \begin{aligned}
 C^+_{1,0;0}(\epsilon)&=6\epsilon^2-8,\\
 C^-_{1,1;1}(\epsilon)&=-8,
 &C^+_{1,1;1}(\epsilon)&=16\epsilon^2-8,\\
 C^-_{2,1;2}(\epsilon)&=6\epsilon^2-4,
 &C^+_{2,1;2}(\epsilon)&=30\epsilon^2-4.
 \end{aligned}
\]
In the range $\epsilon^2\geq\frac12$, the last expression is
at least $11$, while $C^+_{1,1;1}$ is nonnegative and vanishes
exactly at $\epsilon^2=\frac12$. The remaining three expressions
give precisely~\eqref{eq:quaternionic-negative-blocks} and the
other two zero cases in~\eqref{eq:quaternionic-zero-blocks}.

When $B=0$ and $(b,x)=(1,2)$, one has $r=3$, so both vertical
Casimir terms are strictly positive. When $(b,x)=(0,2)$, one
has $r=1$ and $(a,b)=(2,0)$; the plus block is positive and the
minus block is identically zero. Its multiplicity is
$d_{2,0}=10$, proving~\eqref{eq:quaternionic-permanent-zero}.
This also excludes all further zero eigenvalues.
\end{proof}

\begin{theorem}\label{thm:quaternionic-index}
The Morse index and nullity of
$\phi:(\bbS^7,g_\epsilon)\to(\bbS^4,h_1)$, with
$\epsilon^2\geq\frac12$, are given by
\begin{equation}\label{eq:quaternionic-index}
 \mathrm{Ind}(\phi)=
 \begin{cases}
 5,&\epsilon^2\geq\frac43,\\
 13,&\frac23\leq\epsilon^2<\frac43,\\
 45,&\frac12\leq\epsilon^2<\frac23.
 \end{cases}
\end{equation}
The nullity is $10$ away from the transition values, and
\begin{equation}\label{eq:quaternionic-nullity}
 \begin{array}{c|ccc}
 \epsilon^2&\frac43&\frac23&\frac12\\
 \hline
 \mathrm{Null}(\phi)&18&42&25 .
 \end{array}
\end{equation}
\end{theorem}
\begin{proof}
Proposition~\ref{prop:quaternionic-spectrum} identifies all negative
and zero blocks. The dimension formula~\eqref{eq:quaternionic-dimension}
gives
\[
 d_{1,0}=4,\qquad d_{1,1}=5,\qquad
 d_{2,1}=16,\qquad d_{2,0}=10.
\]
Thus the three blocks which can be negative have real multiplicities
\[
 \begin{aligned}
 C^+_{1,0;0}:&\quad 2d_{1,0}=8,\\
 C^-_{1,1;1}:&\quad d_{1,1}=5,\\
 C^-_{2,1;2}:&\quad 2d_{2,1}=32.
 \end{aligned}
\]
Summing the strictly negative contributions gives the index values
$5$, $5+8=13$ and $5+8+32=45$ in the respective ranges
of~\eqref{eq:quaternionic-index}.

The permanent zero block $C^-_{2,0;1}$ has multiplicity
$d_{2,0}=10$. At $\epsilon^2=\frac43$ and
$\epsilon^2=\frac23$, the blocks $C^+_{1,0;0}$ and
$C^-_{2,1;2}$ add $2d_{1,0}=8$ and $2d_{2,1}=32$,
respectively. At $\epsilon^2=\frac12$, the additional zero
block $C^+_{1,1;1}$ has multiplicity $3d_{1,1}=15$.
The corresponding nullities are therefore $10+8=18$,
$10+32=42$ and $10+15=25$, proving~\eqref{eq:quaternionic-nullity}.
\end{proof}

For comparison with the canonical variation parameter, the same
index formula is recorded in Table~\ref{tab:morse-eps-thresholds}.
\begin{table}[htbp]
 \centering
 \begin{tabular}{ccc}
 \toprule
 Range of $\epsilon^2$ & Range of $t=(2\epsilon^2)^{-1}$ & $\mathrm{Ind}(\phi)$\\
 \midrule
 $[\frac43,\infty)$ & $(0,\frac38]$ & $5$\\
 $[\frac23,\frac43)$ & $(\frac38,\frac34]$ & $13$\\
 $[\frac12,\frac23)$ & $(\frac34,1]$ & $45$\\
 \bottomrule
 \end{tabular}
 \caption{Quaternionic index in the fiber-squashing range.}
 \label{tab:morse-eps-thresholds}
\end{table}

\begin{corollary}\label{cor:quaternionic-distinguished-metrics}
For the round metric, the normal $\Sp(2)$-homogeneous metric, and
the non-round Einstein metric, respectively, we have
\[
 \begin{array}{c|ccc}
 &\text{Round}&\text{Normal homogeneous}&\text{Non-round Einstein}\\ \hline
 t&1&\frac12&\frac15\\[1ex]
 \epsilon^2&\frac12&1&\frac52\\ [1ex]
 \mathrm{Ind}(\phi)&45&13&5\\ [1ex]
 \mathrm{Null}(\phi)&25&10&10 .
 \end{array}
\]
\end{corollary}
The Einstein parameter is the usual quaternionic squashing
parameter; see~\cite{Ziller1982}.

\begin{remark}\label{rem:quaternionic-conformal-fields}
The permanent negative block $\pi^G_{(1,1)}$ has real dimension $5$
and eigenvalue $-8$. It is realised by the pullbacks of the gradient
conformal fields on the target sphere. Indeed, if $Y$ is such a
field on a round $\bbS^n(r)$ of radius $r$ and $\phi$ is a harmonic Riemannian
submersion, then
\[
 (\nabla^\phi)^\ast\nabla^\phi(Y\circ\phi)=r^{-2}(Y\circ\phi),
 \qquad
 \sum_i R(Y\circ\phi,d\phi(e_i))d\phi(e_i)
       =(n-1)r^{-2}(Y\circ\phi).
\]
Thus $\cJ_\phi(Y\circ\phi)=(2-n)r^{-2}(Y\circ\phi)$.
For the round domain, the
block associated with $\pi^G_{(1,0)}$ has eigenvalue $-5$ and real
multiplicity $8$, and is
realised by $d\phi(X_v)$, where
$X_v(p)=v-\langle v,p\rangle p$, $v\in\R^8$, on the unit $\bbS^7$,
and the brackets here denote the Euclidean inner product;
see the conformal-variation formula in~\cite{ElSoufi93}.
The remaining $32$ negative directions belong to $\pi^G_{(2,1)}$
and have eigenvalue $-1$. Thus the domain and target conformal
fields alone do not account for the full round index.
For a non-round Berger metric, the fields $X_v$ are not generally
conformal; the eigenvalue of the corresponding representation
block is $6\epsilon^2-8$, not a metric-independent conformal
eigenvalue.
\end{remark}

\begin{remark}\label{rem:quaternionic-index-rigidity}
The minimum index $5$ in Theorem~\ref{thm:quaternionic-index}
attains El Soufi's lower bound for harmonic maps into $\bbS^4$
whose stress-energy tensor is nonnegative everywhere and positive
definite somewhere~\cite[Theorem~2.5]{ElSoufi93}.
Indeed, writing $\mathcal H$ and $\mathcal V$ for the horizontal
and vertical distributions of the Hopf Riemannian submersion, one has
$e(\phi)=\frac12|d\phi|^2=2$ and
\[
 S_\phi=e(\phi)g_\epsilon-\phi^*h_1
       =g_\epsilon\rvert_{\mathcal H}
         +2g_\epsilon\rvert_{\mathcal V}>0.
\]
A constant rescaling of the target does not change the index,
so the unit-sphere convention of the cited theorem gives the
same bound here. Whether equality characterises the Hopf map,
up to domain and target isometries, remains a rigidity question.
Without the stress-energy hypothesis, the general bound is only
$\operatorname{Ind}(\psi)\geq2$; equality forces the image to lie
in a totally geodesic $\bbS^2\subset\bbS^4$
\cite[Theorem~2.4]{ElSoufi93}.
The existence of nonconstant harmonic maps from these Berger
spheres into $\bbS^2$ is a related question not settled by the
present calculation.
\end{remark}

\subsection{The octonionic Hopf fibration}\label{subsec:octonionic-Hopf}
Let $G=\Spin(9)$, $K=\Spin(8)$ and $H=\Spin(7)$.
The embedding $K\subset G$ is the standard one induced by
$\SO(8)\subset\SO(9)$. The embedding $H\subset K$ requires a
different choice: $H$ is the stabiliser of a unit vector in one
of the two real half-spin representations of $K$. We choose this
half-spin representation to have highest weight $\varpi_4^K$.
Then
\[
 G/H\cong\bbS^{15},\qquad K/H\cong\bbS^7,\qquad G/K\cong\bbS^8,
\]
and the projection $\phi:G/H\to G/K$ is the octonionic Hopf
fibration. This choice of $H$ is essential: the vector-stabiliser
embedding $\Spin(7)\subset\Spin(8)$ would not give the required
action on $\bbS^{15}$; see~\cite{OrneaPartonPiccinniVuletescu2013}.

We choose
\[
 \bangle{X,Y}=-\frac12\tr_{\R^9}(XY)
\]
on $\fg$ and fix $\delta=2$. Writing
$g_\epsilon=g_{\epsilon,2}$, we have
\begin{equation}\label{eq:octonionic-metric}
 g_\epsilon=\frac14\bangle{\cdot,\cdot}\rvert_{\fq}
       +\epsilon^{-2}\bangle{\cdot,\cdot}\rvert_{\fp},
 \qquad
 g_{\mathrm{rd}}=\frac14\bangle{\cdot,\cdot}\rvert_{\fq}
       +\bangle{\cdot,\cdot}\rvert_{\fp}.
\end{equation}
Thus the unit round domain corresponds to $\epsilon=1$, the target
$(\bbS^8,h_2)$ has radius $\frac12$, and
\begin{equation}\label{eq:octonionic-t-epsilon}
 t=\epsilon^{-2},\qquad
 0<t\leq1\quad\Longleftrightarrow\quad\epsilon\geq1
\end{equation}
in the convention~\eqref{eq:Hopf-short-fibre-convention}.
We restrict the octonionic example to this fiber-squashing range.
These are the octonionic Hopf--Berger metrics; their homogeneous
normalisations are also discussed in~\cite{Bettiol2022-vw}.
The fibres remain totally geodesic, and $\phi$ is harmonic throughout
this range.

For completeness, the factors in~\eqref{eq:octonionic-metric} can
be checked directly in the real spin representation
$\Sigma_9\cong\R^{16}$. If $p\in\Sigma_9$ is a unit vector with
stabiliser $H$, the orbit map sends $X\in\fm$ to $Xp$.
Let $E_{ij}\in\mathfrak{so}(9)$ have entries $1$ at $(i,j)$,
$-1$ at $(j,i)$ and zero elsewhere, for $i<j$.
The unit element $E_{i9}\in\fq$ acts by $\frac12c_ic_9$, where
$c_1,\ldots,c_9$ are symmetric orthogonal Clifford generators on
$\Sigma_9$ satisfying $c_ic_j+c_jc_i=2\delta_{ij}\id$;
so $|E_{i9}p|^2=\frac14$. The spin Casimir is $9$ with the chosen
inner product. Since $\dim\fq=8$, $\dim\fp=7$, and the
$\fp$-summand is $H$-irreducible, its orbit metric factor $A$
satisfies $8\cdot\frac14+7A=9$, whence $A=1$.

We now fix a common maximal torus of $G$ and $K$ and identify
its weight space with $\R^4$. Tuples denote coefficients in
the orthonormal basis $e_1,\ldots,e_4$.
The positive roots of $G$ are $e_i\pm e_j$ for $i<j$, together
with $e_i$. The positive roots of $K$ are $e_i\pm e_j$ for $i<j$.
Consequently,
\[
 \rho_G=\left(\frac72,\frac52,\frac32,\frac12\right),
 \qquad \rho_K=(3,2,1,0).
\]
The fundamental weights of $K$ used below are
\[
 \varpi_1^K=e_1,\qquad \varpi_2^K=e_1+e_2,\qquad
 \varpi_3^K=\frac12(e_1+e_2+e_3-e_4),\qquad
 \varpi_4^K=\frac12(e_1+e_2+e_3+e_4).
\]
An irreducible $G$-representation $\pi^G_\lambda$ has highest
weight
\[
 \lambda=(\lambda_1,\lambda_2,\lambda_3,\lambda_4),\qquad
 \lambda_1\geq\lambda_2\geq\lambda_3\geq\lambda_4\geq0,
\]
where all coordinates are integers or all are half-integers.
An irreducible $K$-representation $\pi^K_\mu$ has highest weight
$\mu_1\geq\mu_2\geq\mu_3\geq|\mu_4|$, in the corresponding
integer or half-integer weight lattice.
The Casimir constants are
\begin{align}
 c_G(\lambda)&\coloneqq\bangle{\lambda+2\rho_G,\lambda}
  =\lambda_1(\lambda_1+7)+\lambda_2(\lambda_2+5)
    +\lambda_3(\lambda_3+3)+\lambda_4(\lambda_4+1),
       \label{eq:cas-Spin9}\\
 c_K(\mu)&\coloneqq\bangle{\mu+2\rho_K,\mu}
  =\mu_1(\mu_1+6)+\mu_2(\mu_2+4)
    +\mu_3(\mu_3+2)+\mu_4^2.\label{eq:cas-Spin8}
\end{align}
The isotropy representation $\Ad_\fq$ is of real type. Its associated
complex representation $\pi_\fq$ is the vector representation of
highest weight $\mu_\fq=\varpi_1^K$, while the highest weight of
the complexified adjoint representation is $\mu_\fg=(1,1,0,0)$.
\begin{lemma}\label{lem:octonionic-casimir-dimension}
The isotropy contribution to $C_\lambda$ is
\begin{equation}\label{eq:octonionic-Cp}
 c_K(\mu_\fq)-c_G(\mu_\fg)
   =c_K(\varpi_1^K)-c_G(1,1,0,0)=7-14=-7.
\end{equation}
The dimension of the irreducible $G$-representation of highest
weight $\lambda$ is
\begin{equation}\label{eq:octonionic-dimension}
 d_\lambda\coloneqq\dim_\C\cH^G_\lambda=
 \prod_{i=1}^4\frac{\lambda_i+\rho_i}{\rho_i}
 \prod_{1\leq i<j\leq4}
 \frac{(\lambda_i+\rho_i)^2-(\lambda_j+\rho_j)^2}
      {\rho_i^2-\rho_j^2},
 \qquad \rho=(\rho_1,\ldots,\rho_4)=\rho_G.
\end{equation}
\end{lemma}
\begin{proof}
The first identity follows from~\eqref{eq:cas-Spin9} and
\eqref{eq:cas-Spin8}; the second is the Weyl dimension formula
for the positive roots of $G$ listed above.
\end{proof}

We now determine the $K$-components contributing to
$\Hom_H(\cH^G_\lambda,\fq_\C)$, starting with the branching rule.
\begin{theorem}[Branching rule, {\cite[Chapter~8]{GoodmanWallach2009}}]
\label{thm:octonionic-branching}
Let $G=\Spin(9)$ and let $K=\Spin(8)$ be the standard subgroup.
For a dominant $G$-weight $\lambda$, the restriction
$\pi^G_\lambda\rvert_K$ is the multiplicity-free sum of the
representations $\pi^K_\mu$ satisfying
\begin{equation}\label{eq:octonionic-branching}
 \lambda_1\geq\mu_1\geq\lambda_2\geq\mu_2
 \geq\lambda_3\geq\mu_3\geq\lambda_4\geq|\mu_4|.
\end{equation}
All coordinates of $\lambda$ and $\mu$ must belong to the same
integer or half-integer weight lattice.
\end{theorem}

\begin{lemma}\label{lem:octonionic-tensor-product}
For $\ell\in\mathbb Z_{\geq0}$, put $s=\ell/2$ and define
\[
 \mu_\ell^+=(s+1,s,s,s),\qquad
 \mu_\ell^-=(s,s,s,s-1),
\]
where the minus type is defined only for $\ell\geq1$.
The weights $\lambda\in\Lambda_G^+$ satisfying
\[
 \Hom_H(\cH^G_\lambda,\fq_\C)\neq0
\]
are exactly the weights in the following four families, where
$n\in\mathbb Z_{\geq0}$:
\begin{equation}\label{eq:octonionic-contributing-weights}
 \begin{array}{c|c|c|c}
 \text{Family}&\lambda&\text{Range}&\text{Intermediate type}\\ \hline
 A^+&(s+1+n,s,s,s)&\ell\geq0&\mu_\ell^+\\
 B^+&(s+1+n,s+1,s,s)&\ell\geq0&\mu_\ell^+\\
 A^-&(s+n,s,s,s)&\ell\geq1&\mu_\ell^-\\
 B^-&(s+n,s,s,s-1)&\ell\geq2&\mu_\ell^-.
 \end{array}
\end{equation}
For a fixed $\lambda$, its relevant $K$-components are precisely
all $\pi^K_{\mu_\ell^\pm}$ arising from its occurrences in
\eqref{eq:octonionic-contributing-weights}. Each occurs once;
the sum $\tau_\lambda$ of all remaining components, on the
representation space $\cH_{\tau_\lambda}$, satisfies
$\Hom_H(\cH_{\tau_\lambda},\fq_\C)=0$.

The highest weights of the $K$-representations with nonzero
$H$-fixed subspace are exactly
\begin{equation}\label{eq:octonionic-spherical}
 \nu_\ell\coloneqq\ell\varpi_4^K,\qquad
 \ell\geq0,\qquad
 (\cH^K_{\nu_\ell})^H\cong\C.
\end{equation}
The corresponding tensor products with $\pi_\fq=\pi^K_{\varpi_1^K}$ satisfy
\begin{equation}\label{eq:octonionic-tensor}
 \pi^K_{\nu_\ell}\otimes\pi_\fq
 \cong\pi^K_{\mu_\ell^+}\oplus\pi^K_{\mu_\ell^-},
 \qquad
 \mu_\ell^+=\varpi_1^K+\ell\varpi_4^K,\qquad
 \mu_\ell^-=\varpi_3^K+(\ell-1)\varpi_4^K,
\end{equation}
with the minus summand omitted when $\ell=0$.
Equivalently, the final types with nonzero $H$-fixed subspace in
$(\pi^K_\mu)^\ast\otimes\pi_\fq$ are precisely the $\pi^K_{\nu_\ell}$
for which $\mu=\mu_\ell^+$ or $\mu=\mu_\ell^-$.
Each occurs once, and the $H$-fixed parts of all other summands vanish.
The relevant Casimir constants are
\begin{equation}\label{eq:octonionic-intermediate}
 \begin{aligned}
 c_K(\mu_\ell^+)&=\ell^2+7\ell+7,&
 c_K(\mu_\ell^-)&=\ell^2+5\ell+1,&
 c_K(\nu_\ell)&=\ell(\ell+6).
 \end{aligned}
\end{equation}
\end{lemma}
\begin{proof}
By the choice of $H$, the space $K/H$ is the unit sphere in the
half-spin representation of highest weight $\varpi_4^K$.
Triality identifies its degree $\ell$ spherical harmonics with
$\pi^K_{\nu_\ell}$. The spherical harmonic decomposition of $L^2(K/H)$
is multiplicity-free and contains exactly these representations;
Frobenius reciprocity proves~\eqref{eq:octonionic-spherical}.

The vector representation $\pi_\fq$ is minuscule (its weights form
a single orbit under the Weyl group of $K$, each with multiplicity
one), with weights $\pm e_i$. Its tensor product rule gives one summand for each
dominant weight obtained by adding one of these weights to
$\ell\varpi_4^K$. The only possibilities are
$\ell\varpi_4^K+e_1$ and, for $\ell\geq1$,
$\ell\varpi_4^K-e_4$. This proves~\eqref{eq:octonionic-tensor}.
All $K$-representations here are self-dual. Tensor reciprocity
therefore gives
\[
 \dim\Hom_K\bigl(\pi^K_{\nu_\ell},(\pi^K_\mu)^\ast\otimes\pi_\fq\bigr)
 =\dim\Hom_K\bigl(\pi^K_\mu,\pi^K_{\nu_\ell}\otimes\pi_\fq\bigr).
\]
Together with~\eqref{eq:octonionic-spherical}, this identifies
all final types with $H$-fixed vectors and shows that only the
intermediate types $\mu_\ell^\pm$ can contribute to
$\Hom_H(\cH^G_\lambda,\fq_\C)$.

Apply Theorem~\ref{thm:octonionic-branching} to these types.
For $\mu_\ell^+$, interlacing forces $\lambda_3=\lambda_4=s$,
$\lambda_2=s$ or $s+1$, and $\lambda_1=s+1+n$.
For $\mu_\ell^-$, it forces $\lambda_2=\lambda_3=s$,
$\lambda_1=s+n$, and $\lambda_4=s$ or $s-1$.
The latter alternative requires $\ell\geq2$, since
$\lambda_4\geq0$. The common weight lattice makes $n$ a
nonnegative integer. These are exactly the four families in
\eqref{eq:octonionic-contributing-weights}, with multiplicity
one by the branching theorem. Finally, substitution in
\eqref{eq:cas-Spin8} gives~\eqref{eq:octonionic-intermediate}.
\end{proof}

Applying Theorem~\ref{thm:J-eigenvalues} and
Lemma~\ref{lem:octonionic-tensor-product}, we obtain the entire
spectrum and its sign classification.
\begin{proposition}\label{prop:octonionic-spectrum}
Assume $\epsilon^2\geq1$. For every triple $(\lambda,\ell,\pm)$
such that $\pi^K_{\mu_\ell^\pm}\subset\pi^G_\lambda\rvert_K$,
the corresponding eigenvalue of the octonionic Hopf map is
\begin{equation}\label{eq:octonionic-spectrum}
 \begin{aligned}
 C^\pm_{\lambda;\ell}(\epsilon)
   &\coloneqq C_\lambda(\mu_\ell^\pm,\nu_\ell)\\
   &=4\big(c_G(\lambda)-c_K(\mu_\ell^\pm)-7\big)
      +\epsilon^2\ell(\ell+6).
 \end{aligned}
\end{equation}
Each block has real multiplicity $d_\lambda$, and these blocks
give the entire spectrum. Multiplicities are added whenever
distinct blocks have the same numerical eigenvalue.

The only blocks which can be negative are
\begin{equation}\label{eq:octonionic-negative-blocks}
 \begin{aligned}
 C^+_{(1,0,0,0);0}(\epsilon)&=-24<0
   &&\text{for every }\epsilon^2\geq1,\\
 C^-_{(\frac12,\frac12,\frac12,\frac12);1}(\epsilon)
   &=7\epsilon^2-20<0
   &&\Longleftrightarrow\quad\epsilon^2<\frac{20}{7},\\
 C^+_{(\frac32,\frac12,\frac12,\frac12);1}(\epsilon)
   &=7\epsilon^2-16<0
   &&\Longleftrightarrow\quad\epsilon^2<\frac{16}{7}.
 \end{aligned}
\end{equation}
All other blocks are nonnegative. Moreover,
\begin{equation}\label{eq:octonionic-permanent-zero}
 C^+_{(1,1,0,0);0}(\epsilon)=0
 \qquad\text{for every }\epsilon^2\geq1,
\end{equation}
with real multiplicity $36$. Apart from this permanent zero block,
zero occurs only when
\begin{equation}\label{eq:octonionic-zero-blocks}
 \begin{aligned}
 C^-_{(\frac12,\frac12,\frac12,\frac12);1}(\epsilon)=0
   &\quad\Longleftrightarrow\quad\epsilon^2=\frac{20}{7},\\
 C^+_{(\frac32,\frac12,\frac12,\frac12);1}(\epsilon)=0
   &\quad\Longleftrightarrow\quad\epsilon^2=\frac{16}{7},\\
 C^-_{(1,1,1,0);2}(\epsilon)=16\epsilon^2-16=0
   &\quad\Longleftrightarrow\quad\epsilon^2=1.
 \end{aligned}
\end{equation}
\end{proposition}
\begin{proof}
Substitute $\delta=2$,~\eqref{eq:octonionic-Cp} and
\eqref{eq:octonionic-intermediate} into
Theorem~\ref{thm:J-eigenvalues}. The branching rule and the
relevant tensor products are multiplicity-free, and each final
$H$-fixed subspace has dimension one by
Lemma~\ref{lem:octonionic-tensor-product}.
This proves~\eqref{eq:octonionic-spectrum} and its multiplicities.
As in the quaternionic case, the complex dimension of each
eigenspace in the complexified bundle equals its real multiplicity.

For the sign analysis, put
$B=c_G(\lambda)-c_K(\mu_\ell^\pm)-7$.
The four families in~\eqref{eq:octonionic-contributing-weights}
give
\begin{equation}\label{eq:octonionic-finite-reduction}
 \begin{aligned}
 A^+:\quad&B=2\ell-6+n(\ell+9)+n^2,&&\ell\geq0,\\
 B^+:\quad&B=3\ell+n(\ell+9)+n^2,&&\ell\geq0,\\
 A^-:\quad&B=3\ell-8+n(\ell+7)+n^2,&&\ell\geq1,\\
 B^-:\quad&B=2\ell-8+n(\ell+7)+n^2,&&\ell\geq2.
 \end{aligned}
\end{equation}
All four expressions are strictly positive for $n\geq1$ in
their respective ranges. For $n=0$ and $\epsilon^2\geq1$,
the eigenvalue is bounded below by $4B+\ell(\ell+6)$.
For $A^+,B^+,A^-,B^-$ these lower bounds are, respectively,
\[
 \ell^2+14\ell-24,\qquad \ell^2+18\ell,\qquad
 \ell^2+18\ell-32,\qquad \ell^2+14\ell-32.
\]
They leave only $A^+$ with $\ell=0,1$, $A^-$ with $\ell=1$,
and $B^-$ with $\ell=2$, together with $B^+$ at $\ell=0$.
The first three give exactly~\eqref{eq:octonionic-negative-blocks}.
The $B^-$ block is $16\epsilon^2-16$, so it is nonnegative
and vanishes only at $\epsilon^2=1$. The $B^+$ block is
identically zero and has multiplicity $d_{(1,1,0,0)}=36$ by
\eqref{eq:octonionic-dimension}. This proves
\eqref{eq:octonionic-permanent-zero} and
\eqref{eq:octonionic-zero-blocks}, and excludes all further zeros.
\end{proof}

\begin{theorem}\label{thm:octonionic-index}
The Morse index and nullity of
$\phi:(\bbS^{15},g_\epsilon)\to(\bbS^8,h_2)$, with
$\epsilon^2\geq1$, are given by
\begin{equation}\label{eq:octonionic-index}
 \mathrm{Ind}(\phi)=
 \begin{cases}
 9,&\epsilon^2\geq\frac{20}{7},\\
 25,&\frac{16}{7}\leq\epsilon^2<\frac{20}{7},\\
 153,&1\leq\epsilon^2<\frac{16}{7}.
 \end{cases}
\end{equation}
The nullity is $36$ away from the transition values, and
\begin{equation}\label{eq:octonionic-nullity}
 \begin{array}{c|ccc}
 \epsilon^2&\frac{20}{7}&\frac{16}{7}&1\\
 \hline
 \mathrm{Null}(\phi)&52&164&120 .
 \end{array}
\end{equation}
\end{theorem}
\begin{proof}
Proposition~\ref{prop:octonionic-spectrum} identifies every
negative and zero block. The dimension
formula~\eqref{eq:octonionic-dimension} gives
\[
 \begin{aligned}
 d_{(1,0,0,0)}&=9,&
 d_{(\frac12,\frac12,\frac12,\frac12)}&=16,&
 d_{(\frac32,\frac12,\frac12,\frac12)}&=128,\\
 d_{(1,1,1,0)}&=84,&d_{(1,1,0,0)}&=36.
 \end{aligned}
\]
These are the real block multiplicities, since the final
$H$-fixed subspaces have dimension one. Summing the strictly
negative contributions gives $9$, $9+16=25$ and
$9+16+128=153$ in the respective ranges
of~\eqref{eq:octonionic-index}.

The permanent zero block has multiplicity $36$.
At $\epsilon^2=\frac{20}{7}$ and $\epsilon^2=\frac{16}{7}$,
the additional zero blocks have multiplicities $16$ and $128$,
respectively. At $\epsilon^2=1$, the additional block has
multiplicity $84$. Thus the exceptional nullities are
$36+16=52$, $36+128=164$ and $36+84=120$, proving
\eqref{eq:octonionic-nullity}.
\end{proof}

For comparison with the canonical variation parameter, the same
index formula is recorded in Table~\ref{tab:octonionic-index-thresholds}.
\begin{table}[htbp]
 \centering
 \begin{tabular}{ccc}
 \toprule
 Range of $\epsilon^2$ & Range of $t=\epsilon^{-2}$ & $\mathrm{Ind}(\phi)$\\
 \midrule
 $[\frac{20}{7},\infty)$ & $(0,\frac7{20}]$ & $9$\\
 $[\frac{16}{7},\frac{20}{7})$ & $(\frac7{20},\frac7{16}]$ & $25$\\
 $[1,\frac{16}{7})$ & $(\frac7{16},1]$ & $153$\\
 \bottomrule
 \end{tabular}
 \caption{Octonionic index in the fiber-squashing range.}
 \label{tab:octonionic-index-thresholds}
\end{table}

\begin{corollary}\label{cor:octonionic-distinguished-metrics}
For the round metric, the normal $\Spin(9)$-homogeneous metric
(up to homothety), and the non-round Einstein metric, respectively,
\[
 \begin{array}{c|ccc}
 &\text{Round}&\text{Normal homogeneous}&\text{Non-round Einstein}\\ \hline
 t&1&\frac14&\frac3{11}\\
 \epsilon^2&1&4&\frac{11}{3}\\
 \mathrm{Ind}(\phi)&153&9&9\\
 \mathrm{Null}(\phi)&120&36&36 .
 \end{array}
\]
\end{corollary}
For the middle column, $g_{\epsilon=2}=\frac14\bangle{\cdot,\cdot}\rvert_\fm$,
which explains the homothety qualification. The Einstein parameter
can also be recovered from the Ricci eigenvalues of the canonical
variation:
\[
 \Ric_{g_t}\rvert_{\fp}=\left(\frac{6}{t}+8t\right)g_t\rvert_{\fp},
 \qquad
 \Ric_{g_t}\rvert_{\fq}=(28-14t)g_t\rvert_{\fq}.
\]
Equating them gives $(t-1)(11t-3)=0$;
see also~\cite{Ziller1982}.

\begin{remark}\label{rem:octonionic-round-fields}
At the round parameter $\epsilon=1$, the negative eigenvalues
are $-24$, $-13$ and $-9$, with multiplicities $9$, $16$ and
$128$, respectively. The first block is realised by target
gradient conformal fields pulled back along $\phi$, and the
second by the sections $d\phi(X_v)$ induced by the conformal
fields $X_v(p)=v-\langle v,p\rangle p$, $v\in\R^{16}$, of
the round unit $\bbS^{15}$, with the Euclidean inner product.
The third block is
$\pi^G_{(\frac32,\frac12,\frac12,\frac12)}$ and is not accounted
for by those two constructions. In particular, counting only
the domain and target conformal directions gives $25$, whereas
the full index is $9+16+128=153$.
The round kernel consists of the blocks of dimensions $36$ and
$84$, in agreement with
$\mathfrak{so}(16)\rvert_{\Spin(9)}
 \cong\mathfrak{spin}(9)\oplus\Lambda^3\R^9$.
\end{remark}

\begin{remark}\label{rem:Hopf-target-rescaling}
The target radii $\frac12$ in these two examples are chosen so
that the Hopf maps are Riemannian submersions at the round
domain parameters. A constant rescaling of the target leaves
its Levi-Civita connection and its curvature tensor of type
$(1,3)$ unchanged, and hence leaves $\cJ_\phi$ unchanged. Therefore
all indices and nullities above are also valid when the target
spheres are given their round metrics of radius $1$.
\end{remark}

\section{The Jacobi index of the Aloff--Wallach projection}
\label{sec:aloff-wallach}

We consider the canonical homogeneous fibration
\begin{equation}\label{eq:AW-fibration}
 \phi_{k,l}:N_{k,l}=\SU(3)/H_{k,l}\longrightarrow
 \SU(3)/K=\C P^2,\qquad xH_{k,l}\longmapsto xK,
\end{equation}
where
\[
\begin{aligned}
 H_{k,l}&=\{\diag(z^k,z^l,z^{-k-l}):z\in\U(1)\},\\
 K&=\mathrm{S}(\U(2)\times\U(1))
   =\{\diag(A,\det(A)^{-1}):A\in\U(2)\}.
\end{aligned}
\]
We assume $\gcd(k,l)=1$ and $(k,l)\ne(0,0)$. This is the standard
fibration of Aloff--Wallach spaces over the symmetric space $\C P^2$;
see~\cite{GudmundssonSvensson2014}. We use
$\bangle{X,Y}=-\frac12\tr(XY)$ on $\mathfrak{su}(3)$ and the
two-parameter metrics $g_{\epsilon,\delta}$ and $h_\delta$ of
Subsection~\ref{ssec:RHF}. Thus, if
$\mathfrak{su}(3)=\mathfrak h_{k,l}\oplus\mathfrak p\oplus\mathfrak q$
is the orthogonal decomposition, then
$\mathfrak p=\mathfrak k\ominus\mathfrak h_{k,l}$ is vertical,
$\mathfrak q=\mathfrak k^\perp$ is horizontal, and
$\dim\mathfrak p=3$, $\dim\mathfrak q=4$.
The projection is a Riemannian submersion with totally geodesic
fibres. Indeed, for $X,Y\in\mathfrak p$ and $Z\in\mathfrak q$,
Nomizu's formula gives
\[
 2\bangle{\nabla_XY,Z}_{\epsilon,\delta}
 =\bangle{[X,Y]_{\mathfrak m},Z}_{\epsilon,\delta}
  -\bangle{[Y,Z]_{\mathfrak m},X}_{\epsilon,\delta}
  +\bangle{[Z,X]_{\mathfrak m},Y}_{\epsilon,\delta}=0,
\]
since $[\mathfrak p,\mathfrak p]\subset\mathfrak k$ and
$[\mathfrak q,\mathfrak p]\subset\mathfrak q$.
Thus the fibres are minimal and the projection is harmonic.

\begin{lemma}\label{lem:AW-Jacobi}
For every $\epsilon,\delta>0$, the Jacobi operator of
$\phi_{k,l}:(N_{k,l},g_{\epsilon,\delta})\to(\C P^2,h_\delta)$,
under the section identification~\eqref{eq:Psi-section-identification},
is
\begin{equation}\label{eq:AW-Jacobi}
 \cJ_{\phi_{k,l}}
 =\delta^2\bigl(D_K-D_G+D_G^\ad-D_K^\ad\bigr)
   -\epsilon^2D_\Theta.
\end{equation}
Equivalently, on $C^\infty_{H_{k,l}}(\SU(3),\mathfrak q)$,
\begin{equation}\label{eq:AW-Jacobi-simplified}
 \cJ_{\phi_{k,l}}
 =\delta^2\bigl(-D_{\mathfrak q}+D_{\mathfrak q}^\ad\bigr)
   -\epsilon^2D_{\Theta,\mathfrak p}.
\end{equation}
Here, for a $\bangle{\cdot,\cdot}$-orthonormal basis
$Y_1,\ldots,Y_4$ of $\mathfrak q$ and
$Y_5,\ldots,Y_7$ of $\mathfrak p$,
\[
 D_{\mathfrak q}=\sum_{a=1}^4Y_a^2,\qquad
 D_{\mathfrak q}^\ad=\sum_{a=1}^4\ad(Y_a)^2,\qquad
 D_{\Theta,\mathfrak p}=\sum_{a=5}^7\Theta(Y_a)^2.
\]
\end{lemma}
\begin{proof}
The base $\SU(3)/K=\C P^2$ is an irreducible symmetric space. The
only use of the geodesic-orbit hypothesis in the proof of
\cite[Theorem~2.1]{LaraComplexHopf} is to make the connection-trace term
$\sum_i\nabla^g_{E_i}E_i$ vanish at the identity coset, where here
$g=g_{\epsilon,\delta}$ and $\{E_i\}$ is a local orthonormal frame.
We verify
that vanishing directly for this fibration, without assuming that
$g_{\epsilon,\delta}$ is geodesic orbit.

For the reductive homogeneous metric, Nomizu's formula gives
$\nabla^g_{X^\ast}X^\ast|_{eH}=U(X,X)$, where
\[
 \bangle{U(X,X),Z}_{\epsilon,\delta}
 =\bangle{[Z,X]_{\mathfrak p\oplus\mathfrak q},X}_{\epsilon,\delta},
 \qquad Z\in\mathfrak p\oplus\mathfrak q.
\]
If $X\in\mathfrak p$, then
$[\mathfrak p,\mathfrak p]\subset\mathfrak k$ and
$[\mathfrak q,\mathfrak p]\subset\mathfrak q$; the first bracket
has zero pairing with $X$ by $\Ad(\SU(3))$-invariance of
$\bangle{\cdot,\cdot}$, and the second is orthogonal to
$\mathfrak p$. If $X\in\mathfrak q$, then
$[\mathfrak q,\mathfrak q]\subset\mathfrak k$ because
$(\SU(3),K)$ is a symmetric pair, while
$[\mathfrak p,\mathfrak q]\subset\mathfrak q$ and its pairing
with $X$ vanishes by the same invariance. Since the metric is a
constant multiple of $\bangle{\cdot,\cdot}$ on each summand,
$U(X,X)=0$ for $X\in\mathfrak p$ or $X\in\mathfrak q$.
Taking an orthonormal frame consisting of vectors in these two
summands therefore makes the entire connection trace vanish.
All other steps in the proof of~\cite[Theorem~2.1]{LaraComplexHopf} use only
the symmetric-space structure of the base and remain unchanged.
This proves~\eqref{eq:AW-Jacobi}.

For the simplified expression, write
$D_G=D_{\mathfrak q}+D_K$ and
$D_G^\ad=D_{\mathfrak q}^\ad+D_K^\ad$.
Moreover, $\Theta(Z)f=0$ for $Z\in\mathfrak h_{k,l}$ and every
$H_{k,l}$-equivariant $f$, by the equivariance condition
\eqref{def:H-eq-functions}. Thus $D_\Theta$ reduces to
$D_{\Theta,\mathfrak p}$, proving~\eqref{eq:AW-Jacobi-simplified}.
\end{proof}

We fix the representation-theoretic notation. For $a,b\in\N_0$,
let $V^G_{a,b}$ be the irreducible $G=\SU(3)$-representation with
highest weight $a\varpi_1^G+b\varpi_2^G$, so $(a,b)$ is its
Dynkin label. These labels differ from the coordinate pairs
used for $\Sp(2)$ in Subsection~\ref{subsec:quaternionic-Hopf}.
On the diagonal torus of $\SU(3)$, write $e_i$ for the coordinate
weights, with $e_1+e_2+e_3=0$; then
$\varpi_1^G=e_1$ and $\varpi_2^G=e_1+e_2$.
For integers $R\geq S$, let $V^K_{R,S}$
be the irreducible $K\cong\U(2)$-representation of highest weight
$(R,S)$: on the diagonal torus $A=\diag(z_1,z_2)$, its highest
weight character is $z_1^Rz_2^S$.
The real isotropy representation is of complex type, with
\[
 \fq_\C=V^K_{2,1}\oplus V^K_{-1,-2}.
\]
These are the two conjugate constituents of the complexified
isotropy representation; both must be included in the spectral count.

\begin{lemma}\label{lem:AW-casimir-dimension}
In the chosen normalization, the Casimir constants are
\begin{align}
 c_G(a,b)&=\frac43\bigl(a^2+ab+b^2+3a+3b\bigr),
    \label{eq:AW-Casimir-SU3}\\
 c_K(R,S)&=\frac43\bigl(R^2-RS+S^2\bigr)+2(R-S).
    \label{eq:AW-Casimir-U2}
\end{align}
The isotropy contribution is
\begin{equation}\label{eq:AW-isotropy}
 c_K(2,1)-c_G(1,1)=6-12=-6,
\end{equation}
and is the same for the conjugate constituent $V^K_{-1,-2}$.
The Weyl dimension formula gives
\begin{equation}\label{eq:AW-dimension}
 d_{a,b}\coloneqq\dim_\C V^G_{a,b}
 =\frac12(a+1)(b+1)(a+b+2).
\end{equation}
\end{lemma}
\begin{proof}
Write a $K$-weight as $R e_1+S e_2$, with $e_3=-e_1-e_2$.
The induced inner product is
\[
 \bangle{R e_1+S e_2,R'e_1+S'e_2}
 =2(RR'+SS')-\frac23(R+S)(R'+S').
\]
The $G$-weight with Dynkin label $(a,b)$ is
$(a+b)e_1+b e_2$, and
$2\rho_G=4e_1+2e_2$, $2\rho_K=e_1-e_2$.
Substitution in $c_L(\mu)=\bangle{\mu,\mu+2\rho_L}$ gives
\eqref{eq:AW-Casimir-SU3}--\eqref{eq:AW-isotropy}.
The dimension formula follows from the three positive roots of $G$.
\end{proof}

We next identify the intermediate $K$-types and their final
$H_{k,l}$-fixed components.
\begin{theorem}[Branching rule, {\cite[Chapter~8]{GoodmanWallach2009}}]
\label{thm:AW-branching}
For $a,b\in\N_0$, the restriction to
$K=\{\diag(A,\det(A)^{-1}):A\in\U(2)\}$ is
\begin{equation}\label{eq:AW-branching}
 V^G_{a,b}\rvert_K
 \cong\bigoplus_{\substack{a+b\geq r\geq b\\ b\geq s\geq0}}
 V^K_{r-c,s-c},\qquad c=a+2b-r-s.
\end{equation}
Every summand occurs once. The actual $K$-weights
$(R,S)=(r-c,s-c)$ can have negative coordinates.
\end{theorem}
\begin{proof}
The $\U(3)$-representation of highest weight $(a+b,b,0)$
restricts under $\U(2)\times\U(1)$ to the multiplicity-free sum
of $V_{r,s}^{\U(2)}\otimes z^c$, with the interlacing
conditions and exponent $c$ in~\eqref{eq:AW-branching}.
On restricting to $K$, the character $z^c$ becomes
$\det(A)^{-c}$, shifting the $\U(2)$-weight to $(r-c,s-c)$.
Distinct pairs $(r,s)$ give distinct shifted weights, proving the claim.
\end{proof}

\begin{lemma}\label{lem:AW-tensor-product}
For each pair of integers $R\geq S$, define the multiset
\begin{equation}\label{eq:AW-tensor-constituents}
 \mathcal N_{R,S}
 =\{(2-S,1-R),\,(-S-1,-R-2)\}
 \,\sqcup\,
 \begin{cases}
 \{(1-S,2-R),\,(-S-2,-R-1)\},&R>S,\\
 \varnothing,&R=S.
 \end{cases}
\end{equation}
Then
\begin{equation}\label{eq:AW-tensor-product}
 (V^K_{R,S})^\ast\otimes\fq_\C
 \cong\bigoplus_{(u,v)\in\mathcal N_{R,S}}V^K_{u,v},
\end{equation}
where repetitions in the multiset are retained. For each final type,
\begin{equation}\label{eq:AW-H-fixed}
 m_{k,l}(u,v)\coloneqq\dim_\C(V^K_{u,v})^{H_{k,l}}
 =\sum_{j=0}^{u-v}\mathbf1_{\{ku+lv+j(l-k)=0\}}.
\end{equation}
In particular,
$\Hom_{H_{k,l}}(V^G_{a,b},\fq_\C)\neq0$ if and only if
some summand $(R,S)$ of~\eqref{eq:AW-branching} has a final
weight $(u,v)\in\mathcal N_{R,S}$ with $m_{k,l}(u,v)>0$.
These are exactly the contributing intermediate and final types;
all remaining final summands have zero $H_{k,l}$-fixed part.
For $H_{1,1}$, the fixed dimensions simplify to
\begin{equation}\label{eq:AW11-H-fixed}
 m_{1,1}(u,v)=
 \begin{cases}
 u-v+1,&u+v=0,\\
 0,&u+v\neq0.
 \end{cases}
\end{equation}
Consequently, only intermediate types with $R+S=3$ or $R+S=-3$
can contribute for $H_{1,1}$.
\end{lemma}
\begin{proof}
Write $E=\C^2$ for the defining representation of $\U(2)$.
Here $\Sym^d E$ is its $d$th symmetric power, and $\det^j$
denotes the one-dimensional character $A\mapsto\det(A)^j$.
For $d=R-S$, one has
\[
 (V^K_{R,S})^\ast=\det^{-R}\otimes\Sym^d E,\qquad
 V^K_{2,1}=\det\otimes E,\qquad
 V^K_{-1,-2}=\det^{-2}\otimes E.
\]
The Clebsch--Gordan formula~\cite[p.~339]{GoodmanWallach2009},
with the determinant characters included, reads
\[
 \Sym^d E\otimes E
 \cong\Sym^{d+1}E\oplus(\det\otimes\Sym^{d-1}E),
\]
where the second summand is omitted for $d=0$.
Tensoring this identity by $\det^{1-R}$ gives the weights
$(2-S,1-R)$ and $(1-S,2-R)$.
Tensoring by $\det^{-R-2}$ gives $(-S-1,-R-2)$ and
$(-S-2,-R-1)$. This proves~\eqref{eq:AW-tensor-product}.

The weights of $V^K_{u,v}$ are $(u-j,v+j)$, once each, for
$0\leq j\leq u-v$. On $H_{k,l}$, their circle weights are
$ku+lv+j(l-k)$. Counting the zero weights proves
\eqref{eq:AW-H-fixed}; combining it with the branching theorem
and the identification
\[
 \Hom_{H_{k,l}}(V^G_{a,b},\fq_\C)
 \cong\bigl((V^G_{a,b})^\ast\otimes\fq_\C\bigr)^{H_{k,l}}
\]
gives the stated selection of components.
For $k=l=1$, every weight has circle weight $u+v$,
proving~\eqref{eq:AW11-H-fixed}.
In~\eqref{eq:AW-tensor-constituents}, the two channels from
$V^K_{2,1}$ have $u+v=3-R-S$, and the two from $V^K_{-1,-2}$
have $u+v=-3-R-S$. This proves the last assertion.
\end{proof}

The operator formula~\eqref{eq:AW-Jacobi} and
Lemma~\ref{lem:AW-tensor-product} now give the spectrum and a
finite classification of all blocks which can be nonpositive.
\begin{proposition}\label{prop:AW-spectrum}
For $\epsilon,\delta>0$, the eigenvalues of $\cJ_{\phi_{k,l}}$ are
\begin{equation}\label{eq:AW-eigenvalues}
 \begin{aligned}
 \Lambda_{a,b;R,S;u,v}
 &=\delta^2 A_{a,b;R,S}+\epsilon^2c_K(u,v),\\
 A_{a,b;R,S}&\coloneqq c_G(a,b)-c_K(R,S)-6,
 \end{aligned}
\end{equation}
where $(R,S)=(r-c,s-c)$ ranges over~\eqref{eq:AW-branching},
$(u,v)\in\mathcal N_{R,S}$, and $m_{k,l}(u,v)>0$.
Each block has real multiplicity $d_{a,b}m_{k,l}(u,v)$.
Coincident eigenvalues have the sum of their block multiplicities.

Put $\eta=\epsilon^2/\delta^2$.
The only intermediate blocks which can give negative eigenvalues
are the following:
\begin{equation}\label{eq:AW-negative-intermediates}
 \begin{array}{c|c|r}
 (a,b)&(R,S)&A_{a,b;R,S}\\ \hline
 (0,0)&(0,0)&-6\\
 (1,0)&(1,0)&-4\\
 (0,1)&(0,-1)&-4\\
 (1,0)&(-1,-1)&-2\\
 (0,1)&(1,1)&-2\\
 (2,0)&(2,0)&-2\\
 (0,2)&(0,-2)&-2\\
 (1,1)&(1,-1)&-2.
 \end{array}
\end{equation}
Let $\mathcal T$ denote these eight quadruples $(a,b,R,S)$.
For each row and each final type with $m_{k,l}(u,v)>0$, one
has $c_K(u,v)>0$ and
\begin{equation}\label{eq:AW-negative-blocks}
 \Lambda_{a,b;R,S;u,v}<0
 \quad\Longleftrightarrow\quad
 \eta<\eta^\ast_{a,b;R,S;u,v},\qquad
 \eta^\ast_{a,b;R,S;u,v}\coloneqq\frac{-A_{a,b;R,S}}{c_K(u,v)}.
\end{equation}
That block is zero exactly when $\eta=\eta^\ast_{a,b;R,S;u,v}$.
There are also two permanent zero blocks, namely
\begin{equation}\label{eq:AW-permanent-zero}
 (a,b)=(1,1),\qquad (R,S)=(2,1)\text{ or }(-1,-2),
 \qquad (u,v)=(0,0).
\end{equation}
They have total real multiplicity $16$ for every $(k,l)$.
All other blocks are strictly positive. For $(k,l)=(1,1)$,
none of the rows in~\eqref{eq:AW-negative-intermediates}
has a final type with nonzero $H_{1,1}$-fixed part.
\end{proposition}
\begin{proof}
Apply Peter--Weyl to the complexified section space, followed by
Theorem~\ref{thm:AW-branching} and
Lemma~\ref{lem:AW-tensor-product}.
On each block of~\eqref{eq:AW-Jacobi}, the Casimir operators
act by the constants in Lemma~\ref{lem:AW-casimir-dimension}.
The isotropy contribution is $-6$, which proves
\eqref{eq:AW-eigenvalues}.
The branching rule is multiplicity-free, and the final fixed
subspace has dimension $m_{k,l}(u,v)$.
The complex dimension of each eigenspace in the complexified
real bundle equals its real multiplicity. Both conjugate
isotropy channels have already been included in
$\mathcal N_{R,S}$, so no additional factor is needed.

For the sign analysis, put $L=a+b$.
As a function of the interlacing variables $(r,s)$,
$c_K(r-c,s-c)$ is convex on the rectangle
$b\leq r\leq a+b$, $0\leq s\leq b$.
Its maximum occurs at a corner. At the four corners the
values of $c_G(a,b)-c_K(r-c,s-c)$ are
\begin{equation}\label{eq:AW-finite-reduction}
 \begin{array}{c|c}
 (r,s)&c_G(a,b)-c_K(r-c,s-c)\\ \hline
 (a+b,b)&2a+4b\\
 (a+b,0)&2(a+b)\\
 (b,b)&4ab+4(a+b)\\
 (b,0)&4a+2b.
 \end{array}
\end{equation}
They are all at least $2L$. Consequently,
\[
 A_{a,b;R,S}\geq2L-6,
\]
so every block with $L\geq4$ is strictly positive.
For $L\leq3$, substitution in~\eqref{eq:AW-branching}
gives exactly the eight negative horizontal contributions
listed in~\eqref{eq:AW-negative-intermediates}; all remaining
horizontal contributions are nonnegative.
The final types in those eight rows have positive Casimir,
which proves~\eqref{eq:AW-negative-blocks} and the corresponding
zero thresholds.

It remains to classify zeros when $A_{a,b;R,S}\geq0$.
The final Casimir vanishes only for $(u,v)=(0,0)$.
By~\eqref{eq:AW-tensor-constituents}, this final type occurs
only when $(R,S)=(2,1)$ or $(-1,-2)$.
In either case $A_{a,b;R,S}=c_G(a,b)-12$, which is zero
only for $(a,b)=(1,1)$.
Indeed, $c_G(a,b)=L^2+4L+(a-b)^2/3$, and direct evaluation
for $L\leq2$ leaves exactly that weight; $L\geq3$ gives
$c_G(a,b)\geq21$.
The adjoint restriction is
\[
 V^G_{1,1}\rvert_K
 \cong V^K_{1,-1}\oplus V^K_{0,0}
       \oplus V^K_{2,1}\oplus V^K_{-1,-2}.
\]
Each conjugate isotropy channel gives one scalar final type,
whose fixed dimension is $m_{k,l}(0,0)=1$.
Since $d_{1,1}=8$, these two blocks have total multiplicity
$2d_{1,1}=16$, proving~\eqref{eq:AW-permanent-zero}.
Finally, none of the intermediate weights in
\eqref{eq:AW-negative-intermediates} has $R+S=\pm3$.
Lemma~\ref{lem:AW-tensor-product} therefore excludes all their
final types for $H_{1,1}$.
\end{proof}

\begin{samepage}
\begin{theorem}\label{thm:AW-index}
For every coprime pair $(k,l)\neq(0,0)$ and every
$\epsilon,\delta>0$, put $\eta=\epsilon^2/\delta^2$.
With $\mathcal T$ and $\eta^\ast_{a,b;R,S;u,v}$ as in
Proposition~\ref{prop:AW-spectrum}, the Morse index and nullity are
\begin{equation}\label{eq:AW-general-index}
 \operatorname{Ind}(\phi_{k,l})
 =\sum_{(a,b,R,S)\in\mathcal T}d_{a,b}
   \sum_{(u,v)\in\mathcal N_{R,S}}
   m_{k,l}(u,v)\,
   \mathbf1_{\{\eta<\eta^\ast_{a,b;R,S;u,v}\}},
\end{equation}
\begin{equation}\label{eq:AW-general-nullity}
 \operatorname{Null}(\phi_{k,l})
 =16+\sum_{(a,b,R,S)\in\mathcal T}d_{a,b}
   \sum_{(u,v)\in\mathcal N_{R,S}}
   m_{k,l}(u,v)\,
   \mathbf1_{\{\eta=\eta^\ast_{a,b;R,S;u,v}\}}.
\end{equation}
Both sums are finite, over the eight rows
of~\eqref{eq:AW-negative-intermediates} and their final types.
\end{theorem}
\end{samepage}
\begin{proof}
Proposition~\ref{prop:AW-spectrum} identifies all negative and
zero blocks and their multiplicities.
The dimension formula~\eqref{eq:AW-dimension} gives
\[
 d_{0,0}=1,\qquad d_{1,0}=d_{0,1}=3,\qquad
 d_{2,0}=d_{0,2}=6,\qquad d_{1,1}=8.
\]
Multiplying by $m_{k,l}(u,v)$ and summing the negative blocks
gives~\eqref{eq:AW-general-index}.
At equality, the same multiplicities add to the permanent
kernel of dimension $2d_{1,1}=2\cdot8=16$, giving
\eqref{eq:AW-general-nullity}.
\end{proof}

\begin{corollary}\label{cor:AW11-index}
For all $\epsilon,\delta>0$, the Aloff--Wallach projection
$\phi_{1,1}:(N_{1,1},g_{\epsilon,\delta})\to(\C P^2,h_\delta)$
is stable, and
\begin{equation}\label{eq:AW11-index}
 \operatorname{Ind}(\phi_{1,1})=0,\qquad
 \operatorname{Null}(\phi_{1,1})=16.
\end{equation}
\end{corollary}
\begin{proof}
Every fixed dimension in the finite sums of
Theorem~\ref{thm:AW-index} is zero for $(k,l)=(1,1)$ by
Proposition~\ref{prop:AW-spectrum}.
Only the two permanent zero blocks remain, each of real
multiplicity $d_{1,1}=8$.
\end{proof}

\begin{remark}[Konishi fibration]\label{rem:AW11-geometric-interpretations}
The projection $\phi_{1,1}$ is the Konishi fibration of $\C P^2$.
Its total space admits a homogeneous $3$-Sasakian metric, and the
fibres are the orbits of the three Reeb vector fields;
see~\cite[Section~5.2]{BallMadnick2022Squashed} and
\cite[Remark~5.4]{AslanKarigiannisMadnick2025}.
Since $H_{1,1}$ is the centre of $K\cong\U(2)$, this is the
principal $\mathrm{SO}(3)$-bundle
\[
 \mathrm{SO}(3)\longrightarrow N_{1,1}
 \xrightarrow{\,\phi_{1,1}\,}\C P^2,
 \qquad K/H_{1,1}\cong\mathrm{SO}(3)\cong\mathbb{R}P^3.
\]
It can also be viewed as the bundle of unitary tangent frames modulo
a common phase, whose fibre is
$\U(2)/\U(1)_{\mathrm{central}}=\mathrm{PU}(2)\cong\mathrm{SO}(3)$.
The same total manifold has the presentation
$(\SU(3)\times\mathrm{SO}(3))/\U(2)$ underlying Wilking's positively
curved normal homogeneous space;
see~\cite[Remark~1.1]{DominguezGonzalezMouille2023}.
Thus~\eqref{eq:AW11-index} gives stability and nullity $16$ for
the Konishi projection throughout the metric family considered here.
\end{remark}

\begin{corollary}\label{cor:AW1m2-index}
Let $\epsilon,\delta>0$ and put $\eta=\epsilon^2/\delta^2$.
For $\phi_{1,-2}:(N_{1,-2},g_{\epsilon,\delta})\to(\C P^2,h_\delta)$,
one has
\begin{equation}\label{eq:AW1m2-index}
 \operatorname{Ind}(\phi_{1,-2})=
 \begin{cases}
 34,&0<\eta<\frac19,\\
 18,&\frac19\leq\eta<\frac13,\\
 2,&\frac13\leq\eta<1,\\
 0,&\eta\geq1,
 \end{cases}
 \qquad
 \operatorname{Null}(\phi_{1,-2})=
 \begin{cases}
 32,&\eta\in\{\frac19,\frac13\},\\
 18,&\eta=1,\\
 16,&\text{otherwise}.
 \end{cases}
\end{equation}
In particular, $\phi_{1,-2}$ is stable if and only if
$\epsilon^2\geq\delta^2$.
\end{corollary}
\begin{proof}
For $(k,l)=(1,-2)$, formula~\eqref{eq:AW-H-fixed} becomes
\[
 m_{1,-2}(u,v)
 =\sum_{j=0}^{u-v}\mathbf1_{\{u-2v-3j=0\}}.
\]
Among the eight intermediate blocks in
\eqref{eq:AW-negative-intermediates}, only $(a,b;R,S)=(0,0;0,0)$
and $(1,1;1,-1)$ have final types with nonzero fixed part.
Their contributing final types and eigenvalues are
\[
 \begin{array}{c|c|c|c}
 (a,b;R,S)&(u,v)&\text{Eigenvalue}&\text{Real multiplicity}\\ \hline
 (0,0;0,0)&(2,1),\,(-1,-2)&6\epsilon^2-6\delta^2&2\\
 (1,1;1,-1)&(2,1),\,(-1,-2)&6\epsilon^2-2\delta^2&16\\
 (1,1;1,-1)&(3,0),\,(0,-3)&18\epsilon^2-2\delta^2&16.
 \end{array}
\]
Each final type has $m_{1,-2}(u,v)=1$; the last column sums the
two final types in each row, using $d_{0,0}=1$ and $d_{1,1}=8$.
The zero thresholds are respectively $\eta=1$, $\frac13$ and
$\frac19$. Counting the negative rows, and adding the zero rows
to the permanent nullity $16$, proves~\eqref{eq:AW1m2-index}.
\end{proof}

\begin{remark}[Unit tangent sphere bundle]\label{rem:AW1m2-unit-tangent}
The projection $\phi_{1,-2}$ is the natural unit tangent bundle
projection
\[
 \bbS^3\longrightarrow T^1\C P^2\longrightarrow\C P^2,
 \qquad (x,v)\longmapsto x;
\]
see~\cite[Section~1.2]{DominguezGonzalezMouille2023}.
Indeed, identifying the tangent space at the origin of $\SU(3)/K$
with $\C^2$, the isotropy action is $v\mapsto\det(A)Av$.
The stabilizer of the unit vector $(1,0)$ consists of
$A=\diag(z,z^{-2})$, whose image in $\SU(3)$ is
$\diag(z,z^{-2},z)$, precisely $H_{1,-2}$.
Consequently, $N_{1,-2}\cong T^1\C P^2$, with
$K/H_{1,-2}\cong\bbS^3$, and~\eqref{eq:AW1m2-index} computes the
index and nullity of this projection for the homogeneous metric
family $g_{\epsilon,\delta}$ specified above.
Although permuting the second and third coordinates conjugates
$H_{1,-2}$ to $H_{1,1}$ and identifies their total spaces, it does
not preserve $K$. The unit tangent bundle projection and the
Konishi projection are therefore distinct: their fibres are
$\bbS^3$ and $\mathbb{R}P^3$, respectively.
\end{remark}

\section{Morse index of the Clifford hypersurface in the quaternionic Hopf--Berger sphere}\label{sec:Clifford-hypersurface}

We compute the Morse index of the Clifford hypersurface for the
Berger metrics in the fiber-squashing range, where the fibre scale
$\tau$ satisfies $0<\tau\leq1$. In the notation of
Subsection~\ref{subsec:quaternionic-Hopf}, the relation between this
scale and the metric parameter is
\begin{equation}\label{eq:quaternionic-fibre-scale}
 \tau^2=t=\frac{1}{2\epsilon^2}.
\end{equation}
Thus the range under consideration is equivalent to
$\epsilon^2\geq\frac12$.

We retain the notation of Subsection~\ref{subsec:quaternionic-Hopf}.
Thus
\[
 G=\Sp(2),\qquad
 H=\{\diag(h,1):h\in\Sp(1)\},\qquad
 K=\Sp(1)\times\Sp(1),
\]
and $G/H\cong\bbS^7$ is equipped with the metric
\[
 g_\epsilon
 =\bangle{\cdot,\cdot}\rvert_{\fq}
  +\epsilon^{-2}\bangle{\cdot,\cdot}\rvert_{\fp}.
\]
In particular, the notation $\fp\oplus\fq$ and the parameter
$\epsilon$ agree with those of
Subsection~\ref{subsec:quaternionic-Hopf}.

Put
\[
 x_0=\frac{1}{\sqrt2}
 \begin{pmatrix}1&1\\1&-1\end{pmatrix}\in G
\]
and consider the $K$-orbit
\begin{equation}\label{eq:Clifford-embedding}
 \Phi:K\longrightarrow G/H,\qquad
 \Phi(k)=kx_0H.
\end{equation}
The stabiliser of $x_0H$ in $K$ is trivial, so $\Phi$ is an
embedding. Its image is the Clifford hypersurface
$\bbS^3(1/\sqrt2)\times\bbS^3(1/\sqrt2)\subset\bbS^7$
in the unit round realisation. The domain
$K\cong\bbS^3\times\bbS^3$ is equipped with the induced metric
defined below.

For $\alpha\in\{\mathbf i,\mathbf j,\mathbf k\}\subset\Im\H$, define
\[
 E_\alpha=(\alpha,\alpha),\qquad
 F_\alpha=(\alpha,-\alpha)
 \quad\text{in }\fk=\fsp(1)\oplus\fsp(1).
\]
After translating $T_{x_0H}(G/H)$ to $T_{eH}(G/H)$, the
corresponding tangent vectors are
\[
 X_\alpha=
 \begin{pmatrix}0&0\\0&\alpha\end{pmatrix}\in\fp,
 \qquad
 Y_\alpha=
 \begin{pmatrix}0&\alpha\\\alpha&0\end{pmatrix}\in\fq.
\]
Choose
\[
 Y=\begin{pmatrix}0&1\\-1&0\end{pmatrix}\in\fq
\]
with the sign for which the formula below holds. Then $Y$ is
orthogonal to the six vectors $X_\alpha,Y_\alpha$, and
\[
 |X_\alpha|_{g_\epsilon}^2=\frac{1}{2\epsilon^2},
 \qquad
 |Y_\alpha|_{g_\epsilon}^2=|Y|_{g_\epsilon}^2=1.
\]
If $E_\alpha^L,F_\alpha^L$ denote the corresponding left-invariant
vector fields on $K$, the induced metric
$h_\epsilon=\Phi^\ast g_\epsilon$ therefore has the global
orthonormal frame
\begin{equation}\label{eq:Clifford-orthonormal-frame}
 e_\alpha=\sqrt2\,\epsilon E_\alpha^L,\qquad
 f_\alpha=F_\alpha^L.
\end{equation}
The left translates of $Y$ determine a global unit normal $N$
along $\Phi$.

Let $B_\Phi$ denote the scalar second fundamental form with respect
to $N$. To compute it, consider the nearby principal orbits
$\Phi_\theta(k)=kx_\theta H$, where
\[
 x_\theta=
 \begin{pmatrix}\cos\theta&\sin\theta\\
                 \sin\theta&-\cos\theta\end{pmatrix},
 \qquad x_{\pi/4}=x_0.
\]
The vector $\partial_\theta\Phi_\theta$ is a horizontal unit
normal. Writing $h_\theta=\Phi_\theta^\ast g_\epsilon$, the metric
on the orbit satisfies the identities below. Here Lie algebra elements
are identified with left-invariant fields and $\delta_{\alpha\beta}$
is the Kronecker symbol:
\[
 \begin{aligned}
 h_\theta(E_\alpha,E_\beta)&=\tau^2\delta_{\alpha\beta},\\
 h_\theta(E_\alpha,F_\beta)&=-\tau^2\cos(2\theta)\delta_{\alpha\beta},\\
 h_\theta(F_\alpha,F_\beta)
   &=\bigl(1+(\tau^2-1)\cos^2(2\theta)\bigr)\delta_{\alpha\beta}.
 \end{aligned}
\]
These identities follow by evaluating the round metric and the
vertical projection on the orbit through
$(\sin\theta,-\cos\theta)\in\H^2$.
Since $B_\Phi=-\frac12\partial_\theta h_\theta|_{\theta=\pi/4}$,
the orthonormal frame~\eqref{eq:Clifford-orthonormal-frame} gives
\begin{equation}\label{eq:Clifford-second-fundamental-form}
 \begin{aligned}
 B_\Phi(e_\alpha,e_\beta)&=0,\\
 B_\Phi(f_\alpha,f_\beta)&=0,\\
 B_\Phi(e_\alpha,f_\beta)
   &=-\frac{\delta_{\alpha\beta}}{\sqrt2\,\epsilon}.
 \end{aligned}
\end{equation}
Consequently, $\Phi$ is minimal and
\begin{equation}\label{eq:Clifford-A-norm}
 |B_\Phi|^2=\frac{3}{\epsilon^2}.
\end{equation}
The Ricci tensor of the ambient Hopf--Berger metric satisfies
\begin{equation}\label{eq:Clifford-normal-Ricci}
 \Ric_{g_\epsilon}(N,N)=12-\frac{3}{\epsilon^2};
\end{equation}
see~\cite{Ziller1982}. Hence, with the nonnegative convention for
$\Delta_{h_\epsilon}$, the Jacobi operator of the two-sided minimal
hypersurface is
\begin{equation}\label{eq:Jacobian-clifford}
 \cJ_\Phi
 =\Delta_{h_\epsilon}
   -|B_\Phi|^2-\Ric_{g_\epsilon}(N,N)
 =\Delta_{h_\epsilon}-12.
\end{equation}
Normal variations $uN$ are identified with real functions $u$ on $K$.
In this section, $\operatorname{Ind}(\Phi)$ and
$\operatorname{Null}(\Phi)$ refer to the index and nullity for area,
computed from this scalar normal Jacobi operator.

We now compute the spectrum of $\Delta_{h_\epsilon}$. For
$k\in\mathbb N_0$, let
$\pi_k:\Sp(1)\to\GL(\cH_k)$ be the irreducible complex
representation of highest weight $k$; thus
$\dim_\C\cH_k=k+1$. Every irreducible complex representation of
$K$ is of the form
\[
 \pi^K_{(k,\ell)}=\pi_k\boxtimes\pi_\ell
 \quad\text{on}\quad
 \cH^K_{(k,\ell)}=\cH_k\otimes\cH_\ell.
\]
Peter--Weyl theory gives
\begin{equation}\label{eq:Clifford-Peter-Weyl}
 L^2(K,\C)
 \cong\widehat\bigoplus_{k,\ell\geq0}
 \cH^K_{(k,\ell)}\otimes
       \big(\cH^K_{(k,\ell)}\big)^\ast.
\end{equation}
The corresponding matrix coefficients are
\[
 f_{v,P}(x)=P\big(\pi^K_{(k,\ell)}(x)^{-1}v\big).
\]
Here $v\in\cH^K_{(k,\ell)}$ and
$P\in(\cH^K_{(k,\ell)})^\ast$.

Put $A_\alpha=(\alpha,0)$ and $B_\alpha=(0,\alpha)$. For the
matrix coefficients above, left-invariant derivatives act on the
dual Peter--Weyl factor. Each $\Sp(1)$-representation is self-dual,
so we identify that factor with $\cH_k\otimes\cH_\ell$ when
computing its spectrum. By~\eqref{eq:Clifford-orthonormal-frame},
the positive Laplacian acts on this factor through
\begin{equation}\label{eq:Clifford-Laplacian-Casimir}
 \begin{aligned}
 L_{\epsilon;k,\ell}
 &=-\sum_{\alpha\in\{\mathbf i,\mathbf j,\mathbf k\}}
 \left[
  2\epsilon^2\,d\pi^K_{(k,\ell)}(A_\alpha+B_\alpha)^2
  +d\pi^K_{(k,\ell)}(A_\alpha-B_\alpha)^2
 \right]\\
 &=(2\epsilon^2+1)
   \big(k(k+2)+\ell(\ell+2)\big)\id
 +(2-4\epsilon^2)\Theta_{k,\ell},
 \end{aligned}
\end{equation}
where
\[
 \Theta_{k,\ell}
 =\sum_{\alpha\in\{\mathbf i,\mathbf j,\mathbf k\}}
   d\pi^K_{(k,\ell)}(A_\alpha)
   d\pi^K_{(k,\ell)}(B_\alpha).
\]
Here the quaternionic generators satisfy
$\sum_\alpha d\pi_k(\alpha)^2=-k(k+2)\id$.
Each generator has squared norm $\frac12$ for
$\bangle{\cdot,\cdot}$, so the orthonormal Casimir used in
Subsection~\ref{subsec:quaternionic-Hopf} is twice this sum.

On the diagonal subgroup $\Sp(1)_\Delta=\{(q,q):q\in\Sp(1)\}\subset K$,
the external product $\pi_k\boxtimes\pi_\ell$ restricts to the ordinary
tensor product $\pi_k\otimes\pi_\ell$. The Clebsch--Gordan formula gives
\begin{equation}\label{eq:Clifford-Clebsch-Gordan}
 \cH_k\otimes\cH_\ell
 =\bigoplus_{r=0}^{\min\{k,\ell\}}\cH_{m_r},
 \qquad m_r=k+\ell-2r.
\end{equation}
On the summand $\cH_{m_r}$, the diagonal Casimir identity
\[
 \sum_\alpha
 d\pi^K_{(k,\ell)}(A_\alpha+B_\alpha)^2
 =
 \sum_\alpha d\pi^K_{(k,\ell)}(A_\alpha)^2
 +\sum_\alpha d\pi^K_{(k,\ell)}(B_\alpha)^2
 +2\Theta_{k,\ell}
\]
shows that
\begin{equation}\label{eq:Clifford-Theta-eigenvalue}
 \Theta_{k,\ell}
 =\frac12\left(
 k(k+2)+\ell(\ell+2)-m_r(m_r+2)
 \right)\id.
\end{equation}
Combining~\eqref{eq:Clifford-Laplacian-Casimir}
and~\eqref{eq:Clifford-Theta-eigenvalue} yields the complete spectrum.

\begin{proposition}\label{prop:Clifford-spectrum}
For $k,\ell\in\mathbb N_0$ and
$0\leq r\leq\min\{k,\ell\}$, put $m_r=k+\ell-2r$. The
corresponding eigenvalue of $\Delta_{h_\epsilon}$ is
\begin{equation}\label{eq:Clifford-Laplacian-eigenvalue}
 \lambda_{k,\ell,r}(\epsilon)
 =2\big(k(k+2)+\ell(\ell+2)\big)
  +(2\epsilon^2-1)m_r(m_r+2).
\end{equation}
Its real multiplicity is
\begin{equation}\label{eq:Clifford-multiplicity}
 (k+1)(\ell+1)(m_r+1).
\end{equation}
Consequently, the corresponding Jacobi eigenvalue is
\begin{equation}\label{eq:Clifford-Jacobi-eigenvalue}
 j_{k,\ell,r}(\epsilon)
 =\lambda_{k,\ell,r}(\epsilon)-12.
\end{equation}
\end{proposition}
\begin{proof}
Equation~\eqref{eq:Clifford-Laplacian-eigenvalue} follows immediately
from~\eqref{eq:Clifford-Laplacian-Casimir} and
\eqref{eq:Clifford-Theta-eigenvalue}. The diagonal summand
$\cH_{m_r}$ has complex dimension $m_r+1$, while the untouched
Peter--Weyl factor in~\eqref{eq:Clifford-Peter-Weyl} has dimension
$(k+1)(\ell+1)$. This proves~\eqref{eq:Clifford-multiplicity}.
Since the Laplacian is a real operator, this complex dimension is
the real multiplicity of the corresponding eigenspace on
$C^\infty(K,\R)$. Finally,
\eqref{eq:Clifford-Jacobi-eigenvalue} follows from
\eqref{eq:Jacobian-clifford}.
\end{proof}

The formula also determines all nonpositive blocks in the fiber-squashing
range. If $\epsilon^2\geq\frac12$, then
\[
 j_{k,\ell,r}(\epsilon)
 \geq2\big(k(k+2)+\ell(\ell+2)\big)-12.
\]
Thus a nonpositive block must have
$(k,\ell)=(0,0),(1,0),(0,1)$ or $(1,1)$.
The constant block has eigenvalue $-12$ and multiplicity $1$.
The pair $(1,0,0),(0,1,0)$ has eigenvalue
$6\epsilon^2-9$ and total multiplicity $8$.
For $(k,\ell)=(1,1)$, the block with $r=1$ is identically zero
and has multiplicity $4$, while the block with $r=0$ has
eigenvalue $16\epsilon^2-8$ and multiplicity $12$.
Using $\tau^2=(2\epsilon^2)^{-1}$ gives the following result.

\begin{corollary}\label{cor:Clifford-index}
For $0<\tau\leq1$, the Clifford hypersurface
\[
 \Phi:(\bbS^3\times\bbS^3,h_\epsilon)
          \longrightarrow(\bbS^7,g_\epsilon)
\]
has Morse index and nullity
\[
 \operatorname{Ind}(\Phi)=
 \begin{cases}
  9,&\frac13<\tau^2\leq1,\\
  1,&0<\tau^2\leq\frac13,
 \end{cases}
 \qquad
 \operatorname{Null}(\Phi)=
 \begin{cases}
  16,&\tau^2=1,\\
  12,&\tau^2=\frac13,\\
  4,&0<\tau^2<1,\quad\tau^2\neq\frac13.
 \end{cases}
\]
In particular, the index drops at $\tau^2=\frac13$, where the
additional eight-dimensional block enters the kernel.
\end{corollary}
These values also follow from the general quaternionic result
of~\cite[Theorem~1.1]{Siffert2026HopfSquashing}; the squared fibre
parameter used there is $t=\tau^2$ in our notation.

\section{The Morse index of the identity on quaternionic Hopf--Berger spheres}
\label{sec:Berger-identity}
We retain the metrics and normalisation of
Subsection~\ref{subsec:quaternionic-Hopf}. Thus
\[
 g_\epsilon=g_H+t\,g_V,\qquad
 t=\frac1{2\epsilon^2},\qquad
 \epsilon^2\geq\frac12\quad\Longleftrightarrow\quad0<t\leq1,
\]
where $g_H$ and $g_V$ are the horizontal
and vertical parts of the unit round metric.
The identity $\iota_\epsilon:(\bbS^7,g_\epsilon)\to
(\bbS^7,g_\epsilon)$ is harmonic. Its index calculation differs
from that of the Hopf projection: the bundle is now the entire
tangent bundle, and horizontal and vertical fields can couple.
We use $t$ in the block calculation and state the resulting
index formula in terms of $\epsilon$.
Throughout this section, $(a,b)$ again denotes the coordinate
weight $ae_1+be_2$ of $G=\Sp(2)$, and $d_{a,b}$ is the
$\Sp(2)$-dimension in~\eqref{eq:quaternionic-dimension}.

\begin{lemma}\label{lem:identity-quadratic-form}
For a closed Riemannian manifold $(M,g)$, let $\iota=\id_M$.
The Jacobi operator and quadratic form of $\iota$ are
\begin{equation}\label{eq:identity-Jacobi}
 \begin{aligned}
 \cJ_\iota&=\nabla^*\nabla-\Ric^\sharp=\Delta_1-2\Ric^\sharp,\\
 Q_g(Y)&\coloneqq\int_M g(\cJ_\iota Y,Y)\,dv_g\\
 &=\frac12\int_M|\mathcal L_Yg|^2\,dv_g
                  -\int_M(\operatorname{div}_gY)^2\,dv_g.
 \end{aligned}
\end{equation}
Here $\nabla$ is the Levi-Civita connection, $\mathcal L$ is the
Lie derivative, $dv_g$ is the Riemannian volume measure, and
$\Ric^\sharp$ is defined by
$g(\Ric^\sharp X,Y)=\Ric(X,Y)$.
The nonnegative Hodge Laplacian is $\Delta_1=dd^*+d^*d$ on
one-forms, identified with vector fields by $g$.
In particular, $Q_g$ is nonnegative on divergence-free fields,
and its zero vectors within that subspace are exactly the
Killing fields.
For the above Berger metric, the Ricci endomorphism acts by
\begin{equation}\label{eq:identity-Berger-Ricci}
 R_{\mathcal H}=12-6t,\qquad
 R_{\mathcal V}=\frac2t+4t
\end{equation}
on the horizontal and vertical distributions, respectively.
\end{lemma}
\begin{proof}
For the identity, the curvature term in the Jacobi operator is
$\Ric^\sharp$. The Weitzenb\"ock formula
$\Delta_1=\nabla^*\nabla+\Ric^\sharp$ gives the operator identity.
Using a local orthonormal frame and summing repeated indices,
integration by parts gives
\[
 \int_M\nabla_iY_j\nabla_jY_i\,dv_g
 =\int_M(\operatorname{div}_gY)^2\,dv_g
  -\int_M\Ric(Y,Y)\,dv_g.
\]
Expanding $\mathcal L_Yg$ proves~\eqref{eq:identity-Jacobi}.
On divergence-free fields, its right-hand side vanishes precisely
when $\mathcal L_Yg=0$.
Finally, O'Neill's curvature formulas for the quaternionic Hopf
submersion give $2/t$ from the two vertical planes through a
unit vertical vector and $4t$ from its mixed planes.
For a horizontal vector, the horizontal planes contribute
$12-9t$, and the three mixed planes contribute $3t$.
This proves~\eqref{eq:identity-Berger-Ricci}.
\end{proof}

For the rest of the section we abbreviate
$g=g_\epsilon$ and $\iota=\iota_\epsilon$.
Norms $\|\cdot\|$ are $L^2(g)$-norms, while $|\cdot|$ denotes
the pointwise norm. Gradients, divergence, volume measures and
the musical maps $\flat,\sharp$ use $g$ unless another metric
is displayed. We write $\nabla_{\mathcal H}f$ for the horizontal
projection of $\nabla^g f$.

The commuting actions of $G=\Sp(2)$ and the right fibre group
$\Sp(1)$ are isometric. Under their product, the complex scalar
functions decompose without multiplicities as
\begin{equation}\label{eq:identity-scalar-modules}
 L^2(\bbS^7,\C)
 =\widehat\bigoplus_{a\geq b\geq0}\mathcal F_{a,b},\qquad
 \mathcal F_{a,b}\cong\cH^G_{(a,b)}\boxtimes
                        \cH^{\Sp(1)}_{a-b}.
\end{equation}
This follows also from the branching theorem in
Subsection~\ref{subsec:quaternionic-Hopf}: the $H$-fixed part
of $\pi^G_{(a,b)}$ comes only from the intermediate type $(0,a-b)$.
Put
\begin{equation}\label{eq:identity-block-parameters}
 x=a-b,\qquad B=b(b+x+3),\qquad
 v=x(x+2),\qquad h=4(B+x),\qquad
 \lambda=\lambda_{a,b}(t)=h+\frac vt.
\end{equation}
The scalar Laplacian acts on $\mathcal F_{a,b}$ by $\lambda$,
with real multiplicity
\begin{equation}\label{eq:identity-scalar-multiplicity}
 M_{a,b}=(x+1)d_{a,b}.
\end{equation}
Indeed, writing $k=a+b$, one has
$h=k(k+6)-x(x+2)$, and the scalar spectrum is
$k(k+6)+(t^{-1}-1)x(x+2)$; see
\cite[Theorem~3.8]{Bettiol2022-vw}, whose fibre-length parameter
is $\sqrt t$ in our convention.

\begin{lemma}\label{lem:identity-Berger-blocks}
All negative eigenspaces of $\cJ_{\iota_\epsilon}$ lie in the
vector-field isotypical components of the scalar types
$\mathcal F_{a,b}$ with $(a,b)\neq(0,0)$.
For $x=0$, this component has multiplicity one and eigenvalue
\begin{equation}\label{eq:identity-basic-block}
 4b(b+3)-24+12t,\qquad b\geq1.
\end{equation}
For $x\geq1$, its multiplicity matrix is
\begin{equation}\label{eq:identity-Berger-matrix}
 \mathsf M_{a,b}(t)=
 \begin{pmatrix}
 H_-&0&-\alpha\\
 0&H_+&\beta\\
 -\alpha&\beta&D
 \end{pmatrix},
 \qquad
 \begin{aligned}
 H_-&=\lambda+2x-20+12t,\\
 H_+&=\lambda-2x-24+12t,\\
 D&=\lambda-\frac4t,\\
 \alpha^2&=16t(x+2)+\frac{8t(x+2)B}{x+1},\\
 \beta^2&=\frac{8txB}{x+1}.
 \end{aligned}
\end{equation}
The entries $\alpha,\beta$ are the nonnegative square roots
of the displayed quantities.
When $b=0$, the middle row and column are omitted.
Each eigenvalue of this matrix contributes $M_{a,b}$ to the
real multiplicity, with multiplicities added at coincidences.
All vector-field types absent from the scalar decomposition
have nonnegative quadratic form.
\end{lemma}
\begin{proof}
Divergence intertwines the product-group actions on vector fields
and functions. Hence an isotypical component whose type is absent
from~\eqref{eq:identity-scalar-modules} consists of divergence-free
fields. Lemma~\ref{lem:identity-quadratic-form} makes its quadratic
form nonnegative. Since $\cJ_{\iota_\epsilon}$ commutes with the
isometries, these components cannot couple to scalar types.

We identify the multiplicities of a scalar type using the
isotropy of the product-group action. At $(0,1)\in\H^2$, the
stabiliser is
\[
 \{(\diag(u,q),q):u,q\in\Sp(1)\},
\]
where the product acts by $(A,q)\cdot p=Apq^{-1}$.
Its horizontal and vertical representations have highest weights
$(1,1)$ and $(0,2)$, respectively.
The branching rule of Theorem~\ref{thm:quaternionic-branching}
and Clebsch--Gordan therefore leave the horizontal intermediate
types $(1,x-1)$, for $x\geq1$, and $(1,x+1)$, for $b\geq1$.
There is one vertical channel, from $(0,x)$, precisely when
$x\geq1$. Each channel has multiplicity one.
For $x=0$, only the horizontal gradient channel remains.
Its scalar eigenvalue is $4b(b+3)$; since its gradient is horizontal,
$\cJ_\iota\nabla f=\lambda\nabla f-2\Ric^\sharp\nabla f$
gives~\eqref{eq:identity-basic-block}.

We compute the matrix for $x\geq1$.
Set $(u_1,u_2,u_3)=(\mathbf i,\mathbf j,\mathbf k)$ and
$\xi_i(p)=p\,u_i$ on $\bbS^7\subset\H^2$.
These are right Hopf Killing fields of unit round length, with
$[\xi_i,\xi_j]=2\sum_k\varepsilon_{ijk}\xi_k$.
Here $\varepsilon_{ijk}$ is the alternating symbol with
$\varepsilon_{123}=1$.
For a real-valued function $f\in\mathcal F_{a,b}$ of $L^2$-norm one, set
\[
 U=\nabla_{\mathcal H}f,\qquad
 V=\sum_{i=1}^3(\xi_if)\xi_i,\qquad
 \nabla^{g_\epsilon}f=U+t^{-1}V.
\]
Then $\|U\|^2=h$ and $\|V\|^2=tv$.
Write $U_-$ and $U_+$ for the orthogonal components of $U$
in the two horizontal channels.

First compute the horizontal diagonal entries.
Identify horizontal fields with sections of $\phi^{-1}T\bbS^4$,
and write $\mathcal D_i=\nabla^\phi_{\xi_i}$ and
$I_iY=(\nabla_Y^{g_{\mathrm{rd}}}\xi_i)_{\mathcal H}$.
The canonical-variation connection formulas give, for horizontal $Y$,
\[
 \nabla^{g_\epsilon}_{\xi_i}Y=\mathcal D_iY+tI_iY,
 \qquad
 \sum_{j=1}^4\big|(
  \nabla^{g_\epsilon}_{E_j}Y)_{\mathcal V}\big|^2
 =3t|Y|^2,
\]
where $E_1,\ldots,E_4$ is an orthonormal horizontal frame.
Let $P_{\mathcal H}$ be the orthogonal projection onto
horizontal fields. Thus the horizontal compression of $\cJ_\iota$ is
\begin{equation}\label{eq:identity-horizontal-compression}
 P_{\mathcal H}\cJ_\iota\rvert_{\mathcal H}
 =\cJ_\phi+12t\,\id-2\sum_{i=1}^3 I_i\mathcal D_i.
\end{equation}
Here the two connection-square contributions are $3t$ each;
combining them with $R_{\mathcal H}=12-6t$ and the Hopf
curvature term $12$ gives the constant $12t$.
In the homogeneous frame, $\mathcal D_i=R_i-I_i$, where $R_i$
is right-translation differentiation along
$Z_i=\diag(0,u_i)\in\fp$ in the second $\Sp(1)$ factor of $K$.
In this frame,
$[R_i,I_j]=0$ and $\sum_iI_i^2=-3\id$.
On the channel with intermediate weight $(1,r)$ and final
right-group weight $x$, the Casimirs consequently give
\[
 -2\sum_i I_i\mathcal D_i=\bigl(r(r+2)-v-3\bigr)\id.
\]
Proposition~\ref{prop:quaternionic-spectrum} gives the Hopf
eigenvalue $c_G(a,b)-2r(r+2)-18+v/t$.
Using $c_G(a,b)=h+2v$ in
\eqref{eq:identity-horizontal-compression}, we obtain
$\lambda+v-r(r+2)-21+12t$.
The choices $r=x-1,x+1$ give $H_-,H_+$ and show that the
horizontal off-diagonal entry is zero.

Next, the Killing property of each $\xi_i$ gives
\[
 \mathcal L_Vg_\epsilon
 =\sum_i\bigl(d(\xi_if)\otimes\xi_i^\flat
                +\xi_i^\flat\otimes d(\xi_if)\bigr),
 \qquad \operatorname{div}V=-vf,
\]
where $\flat$ uses $g_\epsilon$.
The scalar Laplacian commutes with $\xi_i$.
Integration by parts and their commutator relations give
\[
 \sum_{i,j}\int(\xi_j\xi_if)(\xi_i\xi_jf)\,dv
   -\int\Big(\sum_i\xi_i^2f\Big)^2\,dv=-4v.
\]
Substitution in~\eqref{eq:identity-Jacobi} yields
\begin{equation}\label{eq:identity-vertical-energy}
 Q_g(V)=tv\lambda-4v,
\end{equation}
so the normalised vertical entry is $D=\lambda-4/t$.

Finally, $\Delta_1df=d\Delta f$ implies
\begin{equation}\label{eq:identity-gradient-equation}
 \cJ_\iota\nabla f
 =(\lambda-2R_{\mathcal H})U
   +t^{-1}(\lambda-2R_{\mathcal V})V.
\end{equation}
Put $w_\pm=\|U_\pm\|^2$.
The horizontal rows of this identity give the mixed quadratic
entries $-2t(x+2)w_-$ and $2txw_+$, respectively.
Its vertical row, using~\eqref{eq:identity-vertical-energy}, gives
\[
 w_-+w_+=h,\qquad (x+2)w_--xw_+=4v.
\]
Therefore
\[
 w_-=4x+\frac{2xB}{x+1},\qquad
 w_+=\frac{2(x+2)B}{x+1}.
\]
Normalising $U_-,U_+,V$ gives $\alpha,\beta$ in
\eqref{eq:identity-Berger-matrix}. If $b=0$, then $B=w_+=0$
and the second horizontal channel is absent.
Each matrix acts on the multiplicity space of $\mathcal F_{a,b}$;
its real block multiplicity is therefore~\eqref{eq:identity-scalar-multiplicity}.
\end{proof}

\begin{proposition}\label{prop:identity-Berger-signs}
For $0<t\leq1$, the scalar types contributing to the Morse index
are exactly the following:
\begin{enumerate}
 \item $\mathcal F_{1,0}$ has one negative eigenvalue for every $t$,
 with real multiplicity $8$.
 \item $\mathcal F_{1,1}$ has eigenvalue $12t-8$, of real
 multiplicity $5$, negative exactly when $t<2/3$.
 \item For each integer $b\geq1$, put $B_b=b(b+4)$ and define
 \begin{equation}\label{eq:identity-transition-polynomial}
 \begin{aligned}
 A_B(t)&=12t^2+(4B-14)t+3,\\
 C_B(t)&=12t^2+(4B-22)t+3,\\
 P_B(t)&=A_B(t)C_B(t)\bigl(4(B+1)t-1\bigr)\\
 &\qquad-4t^3\bigl(3(B+4)C_B(t)+B A_B(t)\bigr).
 \end{aligned}
 \end{equation}
 The polynomial $P_{B_b}$ has a unique positive root $t_b$,
 satisfying
 \begin{equation}\label{eq:identity-transition-bounds}
 \frac1{4B_b+4}<t_b<\frac1{4B_b}.
 \end{equation}
 The type $\mathcal F_{b+1,b}$ has one negative eigenvalue
 exactly when $t<t_b$, and one zero eigenvalue at $t=t_b$.
 Its real multiplicity is
 \begin{equation}\label{eq:identity-jump-multiplicity}
 M_{b+1,b}=2d_{b+1,b}
   =\frac43(b+1)(b+2)(b+3).
 \end{equation}
\end{enumerate}
All remaining scalar-type blocks are strictly positive.
\end{proposition}
\begin{proof}
For $x=0$,~\eqref{eq:identity-basic-block} leaves only $b=1$
as a possible negative or zero block. Its multiplicity is
$d_{1,1}=5$.
For $x=1,b=0$, the two-by-two matrix has determinant
\[
 -\frac{68t^2-26t+3}{t^2}<0,
\]
since the numerator has negative discriminant and positive
leading coefficient. Thus exactly one eigenvalue is negative,
with multiplicity $2d_{1,0}=8$.

We exclude all $x\geq2$ next. If $b=0$, then
$H_-\geq12t>0$, $D\geq4x+v-4$ and
$\alpha^2=16t(x+2)$. Its Schur complement satisfies
\[
 D-\frac{\alpha^2}{H_-}
 \geq4x+v-4-\frac43(x+2)>0.
\]
If $b\geq1$, then $B\geq x+4\geq6$.
Since $v/t+12t\geq16$, one has
\[
 H_-\geq H_+\geq2(B+x+2)>0,\qquad D\geq4(B+x).
\]
Also $\alpha^2+\beta^2=16t(B+x+2)$, so
\[
 D-\frac{\alpha^2}{H_-}-\frac{\beta^2}{H_+}
 \geq4(B+x)-8t>0.
\]
Both cases are positive definite.

It remains to consider $x=1,b\geq1$.
Here $B=B_b\geq5$, and
\[
 H_-=4B-14+\frac3t+12t,\qquad
 H_+=4B-22+\frac3t+12t,\qquad
 D=4B+4-\frac1t.
\]
The horizontal entries are positive, since $3/t+12t\geq12$.
Furthermore, $\alpha^2=12t(B+4)$ and $\beta^2=4tB$.
Consequently, the sign of the Schur complement is the sign of
\[
 \det\mathsf M_{b+1,b}(t)=t^{-3}P_B(t).
\]
Expansion gives
\begin{align*}
 P_B(t)={}&-9+(12B+144)t+(80B^2-192B-812)t^2\\
 &+(64B^3-512B^2+800B+1808)t^3\\
 &+(320B^2-1216B-816)t^4+384Bt^5.
\end{align*}
Every nonconstant coefficient is strictly positive for $B\geq5$
(as is also seen by writing $B=5+y$, $y\geq0$).
Thus $P_B$ is strictly increasing on $t>0$, starts at $-9$
and has exactly one positive zero.
At $t=(4B+4)^{-1}$, $D=0$ and the Schur complement is negative.
At $t=(4B)^{-1}$, $D=4$; using $H_-\geq18$, $H_+\geq10$
and $\alpha^2\leq27/5$, $\beta^2=1$ shows that the Schur
complement is positive. This proves~\eqref{eq:identity-transition-bounds}
and the claimed inertia on either side of $t_b$.
Finally,~\eqref{eq:quaternionic-dimension} gives
\eqref{eq:identity-jump-multiplicity}.
\end{proof}

\begin{theorem}\label{thm:identity-Berger-index}
For every $\epsilon^2\geq1/2$, set
$E_b=(2t_b)^{-1}$, where $t_b$ is defined in
Proposition~\ref{prop:identity-Berger-signs}.
The Morse index of the identity on $(\bbS^7,g_\epsilon)$ is
\begin{equation}\label{eq:identity-Berger-index}
 \operatorname{Ind}(\iota_\epsilon)
 =8+5\mathbf1_{\{\epsilon^2>3/4\}}
  +\sum_{b=1}^{\infty}\frac43(b+1)(b+2)(b+3)
       \mathbf1_{\{\epsilon^2>E_b\}}.
\end{equation}
The thresholds are strictly increasing and satisfy
\begin{equation}\label{eq:identity-epsilon-thresholds}
 2b(b+4)<E_b<2b(b+4)+2,
\end{equation}
so the sum is finite for every fixed $\epsilon$.
Equivalently, if $N(\epsilon)=\#\{b\geq1:E_b<\epsilon^2\}$, then
\begin{equation}\label{eq:identity-Berger-binomial}
 \operatorname{Ind}(\iota_\epsilon)
 =8\binom{N(\epsilon)+4}{4}
    +5\mathbf1_{\{\epsilon^2>3/4\}}.
\end{equation}
In particular, the index tends to infinity as the fibres collapse.
\end{theorem}
\begin{proof}
Lemma~\ref{lem:identity-Berger-blocks} excludes all nonscalar
isotypical components from the negative spectrum.
Proposition~\ref{prop:identity-Berger-signs} classifies the
remaining negative blocks and gives their real multiplicities.
Since $t=(2\epsilon^2)^{-1}$, their conditions become
$\epsilon^2>3/4$ and $\epsilon^2>E_b$, proving
\eqref{eq:identity-Berger-index}.
The bounds~\eqref{eq:identity-transition-bounds} give
\eqref{eq:identity-epsilon-thresholds}; these intervals are
successive and disjoint, hence $E_b$ is strictly increasing
and tends to infinity.
Finally,
\[
 \frac43(b+1)(b+2)(b+3)=8\binom{b+3}{3},\qquad
 \sum_{b=1}^N8\binom{b+3}{3}=8\binom{N+4}{4}-8.
\]
This proves~\eqref{eq:identity-Berger-binomial} and the divergence
of the index as $\epsilon\to\infty$, or $t\to0$.
\end{proof}

For orientation, the first thresholds are
$E_1\approx11.965233$, $E_2\approx25.987752$ and
$E_3\approx43.993433$. Thus the first index values are
\[
 \begin{array}{c|c}
 \text{Range of }\epsilon^2&\operatorname{Ind}(\iota_\epsilon)\\ \hline
 [1/2,3/4]&8\\
 (3/4,E_1]&13\\
 (E_1,E_2]&45\\
 (E_2,E_3]&125.
 \end{array}
\]
At a threshold, the crossing block is zero and is not counted
in the Morse index.

\begin{corollary}\label{cor:identity-distinguished-metrics}
The identity map has Morse index $8$, $13$ and $13$ for the round,
normal $\Sp(2)$-homogeneous and Jensen metrics, respectively.
At the Jensen parameter, the negative eigenvalues are
$-28/5$ and $-13/5$, with real multiplicities $5$ and $8$.
\end{corollary}
\begin{proof}
The distinguished parameters are $\epsilon^2=1/2,1,5/2$.
Since $E_1>10$, the general formula gives the three indices.
We also give an independent verification of the spectrum of the
identity map under the Jensen metric.
For $t=1/5$,~\eqref{eq:identity-Berger-Ricci} gives
$\Ric=(54/5)g$. Hence $\cJ_\iota=\Delta_1-(108/5)\id$ preserves
gradient and divergence-free fields.
On gradients of scalar eigenfunctions it acts by $\lambda-108/5$;
on divergence-free fields its quadratic form is nonnegative by
Lemma~\ref{lem:identity-quadratic-form}.
The scalar eigenvalues are
\[
 k(k+6)+4x(x+2),\qquad 0\leq x\leq k,\quad k\equiv x\pmod2.
\]
For $k\geq3$ they are at least $27>108/5$.
For $k=1$, only $x=1$ occurs, giving $\lambda=19$ with
multiplicity $2d_{1,0}=8$.
For $k=2$, the case $x=0$ gives $\lambda=16$ with
multiplicity $d_{1,1}=5$, while $x=2$ gives $\lambda=48$.
These are therefore all negative directions, with eigenvalues
$19-108/5=-13/5$ and $16-108/5=-28/5$, and total index $13$.
\end{proof}
The identity map under the Jensen metric thus has index $13$,
whereas the quaternionic
Hopf projection from the same metric has index $5$ by
Corollary~\ref{cor:quaternionic-distinguished-metrics}.

\section*{AI statement}
ChatGPT Astra 6.1 was used to assist in writing parts of the text and
in carrying out the representation-theoretic calculations in Sections~4
and~6.

\bibliographystyle{alpha}
\bibliography{bib}

@article{Urakawa1987-so,
  title     = "Stability of harmonic maps and eigenvalues of the Laplacian",
  author    = "Urakawa, Hajime",
  journal   = "Trans. Am. Math. Soc.",
  publisher = "JSTOR",
  volume    =  301,
  number    =  2,
  pages     = "557",
  month     =  jun,
  year      =  1987
}

@article{Riviere2023-co,
  title     = "Harmonic maps from {$S^3$} into {$S^2$} with low Morse index",
  author    = "Rivi{\`e}re, Tristan",
  journal   = "J. Differential Geom.",
  publisher = "International Press of Boston",
  volume    =  125,
  number    =  1,
  month     =  sep,
  year      =  2023
}

@article{Bettiol2022-vw,
  title     = "The first eigenvalue of a homogeneous {CROSS}",
  author    = "Bettiol, Renato G and Lauret, Emilio A and Piccione, Paolo",
  journal   = "J. Geom. Anal.",
  publisher = "Springer Science and Business Media LLC",
  volume    =  32,
  number    =  3,
  month     =  mar,
  year      =  2022,
  language  = "en"
}

@article{ElSoufi93,
     author = {El Soufi, Ahmad},
     title = {Applications harmoniques, immersions minimales et transformations conformes de la sph\`ere},
     journal = {Compositio Mathematica},
     pages = {281--298},
     year = {1993},
     publisher = {Kluwer Academic Publishers},
     volume = {85},
     number = {3},
     mrnumber = {1214448},
     zbl = {0777.58010},
     language = {fr},
     url = {https://www.numdam.org/item/CM_1993__85_3_281_0/}
}

@article{Ziller1982,
  title = {Homogeneous einstein metrics on spheres and projective spaces},
  volume = {259},
  ISSN = {1432-1807},
  url = {http://dx.doi.org/10.1007/BF01456947},
  DOI = {10.1007/bf01456947},
  number = {3},
  journal = {Mathematische Annalen},
  publisher = {Springer Science and Business Media LLC},
  author = {Ziller,  W.},
  year = {1982},
  month = sep,
  pages = {351–358}
}

@book{GoodmanWallach2009,
  author    = {Goodman, Roe and Wallach, Nolan R.},
  title     = {Symmetry, Representations, and Invariants},
  series    = {Graduate Texts in Mathematics},
  volume    = {255},
  publisher = {Springer},
  address   = {New York},
  year      = {2009},
  doi       = {10.1007/978-0-387-79852-3},
  url       = {https://doi.org/10.1007/978-0-387-79852-3}
}

@article{OrneaPartonPiccinniVuletescu2013,
  author  = {Ornea, Liviu and Parton, Maurizio and Piccinni, Paolo
             and Vuletescu, Victor},
  title   = {{Spin(9)} geometry of the octonionic {Hopf} fibration},
  journal = {Transformation Groups},
  volume  = {18},
  number  = {3},
  pages   = {845--864},
  year    = {2013},
  doi     = {10.1007/s00031-013-9233-x},
  url     = {https://arxiv.org/abs/1208.0899}
}

@misc{LaraComplexHopf,
  author        = {Lara, Luiz},
  title         = {Harmonic maps from {$\mathbb{S}^{2n+1}$} to
                   {$\mathbb{C}P^n$} with low {Morse} index},
  year          = {2026},
  eprint        = {2610.07815},
  archivePrefix = {arXiv},
  primaryClass  = {math.DG},
  note          = {arXiv:2610.07815 [math.DG]},
  doi           = {10.48550/arXiv.2610.07815},
  url           = {https://arxiv.org/abs/2610.07815}
}

@article{GudmundssonSvensson2014,
  author  = {Gudmundsson, Sigmundur and Svensson, Martin},
  title   = {Harmonic morphisms from homogeneous spaces of positive curvature},
  journal = {Math. Proc. Cambridge Philos. Soc.},
  volume  = {157},
  number  = {2},
  pages   = {321--327},
  year    = {2014},
  doi     = {10.1017/S0305004114000334}
}

@misc{Siffert2026HopfSquashing,
  author        = {Siffert, Anna},
  title         = {{Morse} index and bifurcation of {Clifford} hypersurfaces
                   under {Hopf} squashing},
  year          = {2026},
  eprint        = {2611.13854},
  archivePrefix = {arXiv},
  primaryClass  = {math.DG},
  url           = {https://arxiv.org/abs/2611.13854}
}

@misc{BallMadnick2022Squashed,
  author        = {Ball, Gavin and Madnick, Jesse},
  title         = {Associative Submanifolds of Squashed {$3$}-{Sasakian} Manifolds},
  year          = {2022},
  eprint        = {2208.10622},
  archivePrefix = {arXiv},
  primaryClass  = {math.DG},
  url           = {https://arxiv.org/abs/2208.10622}
}

@article{DominguezGonzalezMouille2023,
  author  = {Dom{\'i}nguez-V{\'a}zquez, Miguel and
             Gonz{\'a}lez-{\'A}lvaro, David and Mouill{\'e}, Lawrence},
  title   = {Infinite families of manifolds of positive
             {$k$th}-intermediate {Ricci} curvature with {$k$} small},
  journal = {Mathematische Annalen},
  volume  = {386},
  pages   = {1979--2014},
  year    = {2023},
  doi     = {10.1007/s00208-022-02420-w},
  url     = {https://doi.org/10.1007/s00208-022-02420-w}
}

@article{AslanKarigiannisMadnick2025,
  author  = {Aslan, Benjamin and Karigiannis, Spiro and Madnick, Jesse},
  title   = {Calibrated geometry in hyperk{\"a}hler cones,
             {$3$}-{Sasakian} manifolds, and twistor spaces},
  journal = {Canadian Journal of Mathematics},
  volume  = {77},
  number  = {4},
  pages   = {1163--1221},
  year    = {2025},
  doi     = {10.4153/S0008414X24000282},
  url     = {https://doi.org/10.4153/S0008414X24000282}
}

@article{ElSoufi1995,
  author  = {El Soufi, Ahmad},
  title   = {Indice de {Morse} des applications harmoniques de la sph{\`e}re},
  journal = {Compositio Math.},
  volume  = {95},
  number  = {3},
  pages   = {343--362},
  year    = {1995}
}

@article{LoubeauOniciuc2007,
  author  = {Loubeau, Eric and Oniciuc, Cezar},
  title   = {On the biharmonic and harmonic indices of the {Hopf} map},
  journal = {Trans. Amer. Math. Soc.},
  volume  = {359},
  number  = {11},
  pages   = {5239--5256},
  year    = {2007},
  doi     = {10.1090/S0002-9947-07-03934-7}
}

@misc{GutweinLanglais2026,
  author        = {Gutwein, Dominik and Langlais, Thibault},
  title         = {Deformations of harmonic maps with conical singularities},
  year          = {2026},
  eprint        = {2609.20588},
  archivePrefix = {arXiv},
  primaryClass  = {math.DG},
  note          = {arXiv:2609.20588 [math.DG]},
  doi           = {10.48550/arXiv.2609.20588}
}

@misc{Surkov2026,
  author        = {Surkov, Egor},
  title         = {Harmonic maps from spheres with lowest possible index and their properties},
  year          = {2026},
  eprint        = {2609.35589},
  archivePrefix = {arXiv},
  primaryClass  = {math.DG},
  note          = {arXiv:2609.35589 [math.DG]},
  doi           = {10.48550/arXiv.2609.35589}
}

@misc{ChenWangZhang2026,
  author        = {Chen, Qun and Wang, Guofang and Zhang, Mingwei},
  title         = {Harmonic maps from {$\mathbb{S}^{2n+1}$} into {$\mathbb{CP}^n$} with least {Morse} index},
  year          = {2026},
  eprint        = {2609.37945},
  archivePrefix = {arXiv},
  primaryClass  = {math.DG},
  note          = {arXiv:2609.37945 [math.DG]},
  doi           = {10.48550/arXiv.2609.37945}
}

@article{Tanno1993Identity,
  author  = {Tanno, Shukichi},
  title   = {The identity map as a harmonic map of a
             {$(4r+3)$}-sphere with deformed metrics},
  journal = {Kodai Mathematical Journal},
  volume  = {16},
  number  = {2},
  pages   = {171--180},
  year    = {1993},
  doi     = {10.2996/kmj/1138039781},
  url     = {https://doi.org/10.2996/kmj/1138039781}
}

\end{document}